\documentclass[12pt]{amsart}
\usepackage{amssymb, amsthm}
\usepackage[all]{xy}

\begin{document}

\title[Finitely-Generated Projective Modules]{Finitely-Generated Projective Modules over the $\theta$-deformed 4-Sphere}
\author{Mira A. Peterka}
\address{Department of Mathematics \\
University of California \\ Berkeley, CA 94720-3840}
\email{peterka@math.berkeley.edu}
\thanks{The research reported here was
supported in part by National Science Foundation grant DMS-0500501.}

\maketitle
\thispagestyle{empty}

\newtheorem{theorem}{Theorem}[section]
\newtheorem{cor}[theorem]{Corollary} 
\newtheorem{prop}[theorem]{Proposition} 
\newtheorem{Item}{Item}[section]
\newtheorem{mainlemma}[theorem]{Main Lemma} 
\newtheorem{lemma}[theorem]{Lemma} 
\newtheorem{sublemma}[theorem]{Sublemma} 

\theoremstyle{definition}   %my change
\newtheorem{definition}[theorem]{Definition} 
\newtheorem{example}[theorem]{Example} 
\newtheorem{notation}[theorem]{Notation} 

\renewcommand{\subjclassname}{%
   \textup{2000} Mathematics Subject Classification}

% Fill in the blanks or delete the following

\begin{abstract}
We investigate the ``$\theta$-deformed spheres" $C(S^{3}_{\theta})$ and $C(S^{4}_{\theta})$, where $\theta$ is any real number. We show that all finitely-generated projective modules over $C(S^{3}_{\theta})$ are free, and that $C(S^{4}_{\theta})$ has the cancellation property. We classify and construct all finitely-generated projective modules over $C(S^{4}_{\theta})$ up to isomorphism. An interesting feature is that if $\theta$ is irrational then there are nontrivial ``rank-1" modules over $C(S^{4}_{\theta})$. In that case, every finitely-generated projective module over $C(S^{4}_{\theta})$ is a sum of a rank-1 module and a free module. If $\theta$ is rational, the situation mirrors that for the commutative case $\theta=0$.
\end{abstract}

\maketitle
\allowdisplaybreaks

Of the numerous quantum analogs of spheres that have been discovered (e.g. \cite{woron}, \cite{podles}, \cite{vaksman}, \cite{reshetik}, \cite{bonechi}, \cite{masuda}), the so-called $\theta$-deformed spheres $C(S^{l}_{\theta})$ of Connes and Landi \cite{isospectral} most closely resemble the noncommutative tori $C(T^{l}_{\theta})$. The noncommutative tori are the most studied examples of  ``noncommutative differential manifolds", and have been found to arise naturally in a number of contexts in mathematics and quantum field theory/string theory, for instance in the physics of D-branes and string dualities (e.g. \cite{NCG}, \cite{compactification}, \cite{seiberg}). One might ask whether the $C(S^{l}_{\theta})$ also have physical significance, and indeed there has been much (mostly algebraic) work (see \cite{landivs}, \cite{lanreina}, \cite{brainvs}, \cite{brain}) on developing classical Yang-Mills theory for $C(S^{4}_{\theta})$. Work (\cite{brainvs}, \cite{brain}) characterizing the associated moduli spaces of instanton solutions using a modified ADHM construction already shows relations to (noncommutative) algebraic geometry along lines similar to those from the ADHM construction (see \cite{yangmills}, \cite{adhm}, \cite{donaldson}) for the Moyal-deformed $\mathbb{R}^{4}$ studied by Nekrasov and Schwarz (\cite{nekrasov}, \cite{kapustin}). One would expect that the isomorphism classes of  ``vector bundles"  (finitely-generated projective modules) over $C(S^{4}_{\theta})$ should correspond to the different topological sectors of instanton solutions for the Yang-Mills theory. However, a complete classification and construction (up to isomorphism), of the finitely-generated  projective (right) $C(S^{4}_{\theta})$-modules has not yet appeared. Such a classification and construction is the main result of the present work. 

We note that since the algebras $C(S^{l}_{\theta})$ are noncommutative, we must make a choice between working with right or left $C(S^{l}_{\theta})$-modules. We work exclusively with right modules, and therefore usually simply write ``module"  instead of ``right module". In the case that an algebra $A$ is commutative, we treat all $A$-modules as right $A$-modules.

We find that, for irrational $\theta$, the semigroup $V(C(S^{4}_{\theta}))$ of isomorphism classes of finitely-generated projective modules over $C(S^{4}_{\theta})$ is isomorphic to $\{0\}\cup (\mathbb{N} \times \mathbb{Z})$. For each $(n,s) \in \mathbb{N} \times \mathbb{Z}$, we explicitly construct a $C(S^{4}_{\theta})$-module $N(n,s) \cong M(n,s) \otimes_{C(S^{2})} C(S^{4}_{\theta})$, where $M(n,s)$ is the module of sections of a rank-$n$ complex vector bundle over $S^{2}$ with Chern number $-s$. This gives one of several senses in which the two parameters $n$ and $s$ can be thought of as ``rank" and ``index" respectively. From the isomorphism $N(n,s) \cong N(1,s)\oplus N(n-1,0)$, we see that every finitely-generated projective $C(S^{4}_{\theta})$-module is either of ``rank-1"  or is isomorphic to a direct sum of a ``rank-1"  module and a free module.

Nontrivial finitely-generated projective $C(S^{4}_{\theta})$-modules have already appeared in the literature. All of these have ``rank" $\geq 2$ in our sense, and so cannot be indecomposable if $\theta$ is irrational. In particular, if $\theta$ is irrational, then the noncommutative ``instanton bundle of charge (instanton number) 1" $eC(S^{4}_{\theta})^{4}$, of Connes and Landi \cite{isospectral}, must be isomorphic to $N(2,-1)$, and so $eC(S^{4}_{\theta})^{4}$ is isomorphic to the direct sum $N(1,-1)\oplus C(S^{4}_{\theta})$. The group $K_{0}(C(S^{4}_{\theta}))\cong \mathbb{Z \times Z}$ is generated, however, by the $K_{0}$-classes of the rank-1 free $C(S^{4}_{\theta})$-module and $eC(S^{4}_{\theta})^{4}$.

The modules $p_{(n)}$ of Landi and van Suijlekom \cite{landivs} are isomorphic to $N(n+1, -(1/6)n(n+1)(n+2))$ for $n\geq 1$. Brain and Landi \cite{brain} took representatives from each class of $V(C(S^{4}))$ and generalized them to obtain finitely-generated projective $C(S^{4}_{\theta})$-modules. Additionally, Brain and Landi described a family of (anti-) self-dual  connections on each of the obtained modules. It is clear from noncommutative index calculations that these modules are in distinct classes in $V(C(S^{4}_{\theta}))$. It was not clear, however, whether there might be any finitely-generated projective $C(S^{4}_{\theta})$-modules not isomorphic to any of the modules that Brain and Landi construct.  Our results show that, if $n \geq 2$, and $s \in \mathbb{Z}$ is arbitrary, then $N(n, s)$ is isomorphic to one of Brain and Landi's modules, so there is a bijection between the semigroup consisting of  the isomorphism classes of  the modules that Brain and Landi construct, and the semigroup consisting of the isomorphism class of the free rank-1 $C(S^{4}_{\theta})$-module together with the isomorphism classes of all the $N(n,s)$ with $n \geq 2$.  But none of Brain and Landi's modules are isomorphic to a module of the form $N(1,s$), for $s \neq 0$. In fact, no nontrivial ``line bundles"  over $C(S^{4}_{\theta})$ have appeared previously in the literature.

If $\theta=p/q$ is rational, then $V(C(S^{4}_{p/q})) \cong V(C(S^{4})) \cong \{0, 1\} \cup ((\mathbb{N} \backslash \{1\}) \times \mathbb{Z}))$, as a consequence of the theory of strong Morita equivalence. Again we explicitly construct representatives $N_{p/q}(n, s)$ for each isomorphism class, where $n \geq 2$ can be thought of as rank and $s \in \mathbb{Z}$ can be thought of as index. As in the classical case $\theta=0$, the modules $N_{p/q}(2, s)$ are indecomposable, and $K_{0}(C(S^{4}_{p/q}))$ is generated by the classes of the rank-1 free $C(S^{4}_{p/q})$-module and Connes and Landi's  $eC(S^{4}_{p/q})^{4} \cong N_{p/q}(2,-1)$.

Our results reflect the fact that, for all real $\theta$ and all $n \geq 2$, the natural map from the group of path-components $\pi_{0}(GL_{n}(C(S^{3}_{\theta})))$ of $GL_{n}(C(S^{3}_{\theta}))$ to $K_{1}(C(S^{3}_{\theta}))\cong \mathbb{Z}$ is an isomorphism. A generator for the group $\pi_{0}(GL_{2}(C(S^{3}_{\theta})))$ is given by the matrix  $\begin{pmatrix} z_{2} & z_{1} \\ -\rho z_{1}^{*} & z_{2}^{*} \end{pmatrix}$, where $\rho=\exp(2\pi i \theta)$, and $z_{1}$ and $z_{2}$ are generators of $C(S^{3}_{\theta})$ given in Definition \ref{oddspheres} below. Note that  in the classical case $\theta=0$, the generators $z_{1}$ and $z_{2}$ are the usual complex coordinate functions for $S^{3}$ viewed as embedded in $\mathbb{C}^{2}$ in the standard way, and the matrix $\begin{pmatrix} z_{2} & z_{1} \\ - z_{1}^{*} & z_{2}^{*} \end{pmatrix}$  is the usual generator of $\pi_{3}(SU(2))$.

 If $\theta$ is rational, then $\pi_{0}(GL_{1}(C(S^{3}_{\theta})))\cong 0$, but if $\theta$ is irrational, then $\pi_{0}(GL_{1}(C(S^{3}_{\theta}))) \cong K_{1}(C(S^{3}_{\theta}))$ under the natural map. If $\theta$ is irrational, then a generator for $\pi_{0}(GL_{1}(C(S^{3}_{\theta})))$ is given explicitly by $X:=\exp(2 \pi it)p+1-p$, where $p$ is a Rieffel projection of trace $\theta \mod 1$.  We form the module $N(1,1)$ by using $X$ as a clutching map. If $\theta$ is irrational, then it follows from our results  that $\begin{pmatrix} X & 0 \\ 0& 1 \end{pmatrix}$ must be homotopic through some path in $GL_{2}(C(S^{3}_{\theta}))$ to $\begin{pmatrix} z_{2} & z_{1} \\ -\rho z_{1}^{*} & z_{2}^{*} \end{pmatrix}$. 

\section*{acknowledgments}

Most of the results presented here appeared in somewhat different form in a doctoral dissertation \cite{thesis} I submitted to the University of California, Berkeley. I thank my dissertation supervisor, Marc A. Rieffel, for carefully reading earlier drafts of this paper, and for making many very helpful comments and suggestions.

\section{The $\theta$-deformed Spheres $C(S^{l}_{\theta})$}

In this section we recall the definition of the $\theta$-deformed spheres along with some basic results concerning them that we will need.
\newline

The $\theta$-deformed $l$-sphere $C(S^{l}_{\theta})$ can be obtained (up to isomorphism) in several ways: 
\newline 
1) As the universal C*-algebra generated by certain spherical relations (see Definitions \ref{oddspheres} and \ref{evenspheres}).
\newline
2) As a certain continuous field of noncommutative tori of appropriate dimensions (see Theorems \ref{oddfield} and \ref{evenfield}).
\newline
3) As an example of Rieffel's general deformation quantization \cite{deformation} by actions of $\mathbb{R}^{m}$. (see Proposition 2.9 of \cite{natsume} and Corollary \ref{deformationquant} below)
\newline
4) As a result of the ``$\theta$-deformation" procedure of Connes and Landi on $S^{l}$ \cite{isospectral}.
\newline
5) As a certain fixed-point algebra \cite{connesdv}.
\newline
6) The relations of Definitions \ref{oddspheres} and \ref{evenspheres} are particular solutions to certain homological equations studied by Connes and Landi \cite{isospectral}, and later Connes and Dubois-Violette \cite{connesdv, moduli}. Thus the universal C*-algebra generated by such a solution is identically one of those given in 1) above.
\newline

 In obtaining the results of this work  we will only need to make use of the basic descriptions 1) and 2).
\newline

The $\theta$-deformed 3-sphere $C(S^{3}_{\theta})$ was first studied by Matsumoto \cite{matsumoto}, who defined it as a noncommutative version of the genus-1 Heegaard decomposition of $S^{3}$ as the union of two solid tori glued together along their boundaries. Natsume and Olsen \cite{natsume} introduced the $\theta$-deformed spheres $C(S^{2m-1}_{\theta})$ of arbitrary odd dimension $2m-1 \geq 3$ in work generalizing the classical results of Coburn \cite{coburn} on Toeplitz operators on odd-dimensional spheres.

\begin{definition}
\label{oddspheres}
Let $(\theta)_{ij}$ be an $m \times m$  real-valued skew-symmetric matrix, and let $\rho_{ij}=\exp(2\pi i {\theta_{ij}})$. The $\theta$-deformed $2m-1$-sphere $C(S^{2m-1}_{\theta})$ is the universal C*-algebra generated by $m$ normal elements $z_{1},..., z_{m}$  satisfying the relations \[ z_{1}z_{1}^{*}+...+ z_{m}z_{m}^{*}=1, \quad z_{i}z_{j}=\rho_{ji}  z_{j}z_{i}.\]
\end{definition}

Natsume and Olsen proved that $C(S^{2m-1}_{\theta})$ is a certain field of noncommutative tori:

\begin{definition}
\label{noncommtori}
Let $(\theta)_{ij}$ be an $m \times m$  real-valued skew-symmetric matrix, and let $\rho_{ij}=\exp(2\pi i {\theta_{ij}})$. The $\theta$-deformed $m$-torus $A_{\theta}=C(T^{m}_{\theta})$ is the universal C*-algebra generated by $m$ unitary elements $U_{1},..., U_{m}$  satisfying the relations $ U_{i}U_{j}=\rho_{ji}  U_{j}U_{i}$.
\end{definition}

Let \[S^{m-1}_{+}=\{(t_{1},..., t_{m}) \in S^{m-1} : t_{k} \geq 0  \text{ for all } k \leq m \}\] denote the first ``octant" of $S^{m-1}$.

\begin{theorem}
\label{oddfield}
(Theorem 2.5 of \cite{natsume}). Let $S(i)=\{(t_{1},..., t_{m}) \in S^{m-1}_{+} : t_{i}=0\}$. Let $A_{\theta}(i)$ be the C*-subalgebra of $A_{\theta}$ generated by the set $\{U_{1},...U_{m}\}$ omitting $U_{i}$. Let $\mathfrak{A}$ be the C*-subalgebra of functions in $C(S^{m-1}_{+})$ which map $S(i)$ into $A_{\theta}(i)$ for all $i$. Then $C(S^{2m-1}_{\theta})\cong \mathfrak{A}$.
\end{theorem}

The generator $z_{i}$ corresponds to the map $(t_{1},...,t_{m}) \mapsto t_{i}U_{i}$ under the isomorphism.

Theorem \ref{oddfield} generalizes the identification of $S^{2m-1}$ with the fibered space over $S^{m-1}_{+}$ that comes from the canonical action of the $m$-torus $T^{m}$ on $S^{2m-1}$.

Of course $S^{m-1}_{+}$ is homeomorphic to the standard $m-1$-simplex $\Delta^{m-1}$ via the map $(t_{1},..., t_{m}) \mapsto (t_{1}^{2},..., t_{m}^{2})$, so we may identify $S^{1}_{+}$ with the unit interval $[0,1]$ by composing the map $S^{1}_{+} \approx \Delta^{1}$ with the map $\Delta^{1} \approx [0,1]: (t_{1}, t_{2}) \rightarrow (t_{1}, 0)$. Doing so, and then composing with the isomorphism of Theorem \ref{oddfield}, we recover Matsumoto's \cite{matsumoto2} earlier characterization of $C(S^{3}_{\theta})$ as the continuous field \[\mathcal{A}=  \{f \in C([0,1], C(T^{2}_{\theta})) : f(0) \in C^{*}(U_{2}), f(1) \in C^{*}(U_{1})\}.\]

Here the generator $z_{1}$ corresponds to the map $t \rightarrow \sqrt{t}U_{1}$, while $z_{2}$ corresponds to the map $t \rightarrow \sqrt{1-t}U_{2}$.
\newline

The study of even-dimensional $\theta$-deformed spheres was initiated by Connes and Landi \cite{isospectral} in work motivated by  considerations of cyclic homology.

\begin{definition}
\label{evenspheres}
Let $(\theta)_{ij}$ be an $m \times m$  real-valued skew-symmetric matrix, and let $\rho_{ij}=\exp(2\pi i {\theta_{ij}})$. The $\theta$-deformed $2m$-sphere $C(S^{2m}_{\theta})$ is the universal C*-algebra generated by $m$ normal elements $z_{1},..., z_{m}$ and a hermetian element $x$  satisfying the relations \[ z_{1}z_{1}^{*}+...+ z_{m}z_{m}^{*}+ x^{2}=1, \quad z_{i}z_{j}=\rho_{ji}  z_{j}z_{i}, \quad [x, z_{i}]=0. \]
\end{definition}

Of course $C(S^{2m}_{\theta})$ is also a continuous field of noncommutative tori, but I have not seen a proof of this in the literature. Towards the proof, it will be useful for us to first consider even-dimensional $\theta$-deformed closed balls:

\begin{definition}
\label{evenball}
Let $(\theta)_{ij}$ be an $m \times m$  real-valued skew-symmetric matrix, and let $\rho_{ij}=\exp(2\pi i {\theta_{ij}})$. The $\theta$-deformed $2m$-ball $C(D^{2m}_{\theta})$ is the universal C*-algebra generated by $m$ normal elements $w_{1},..., w_{m}$ and a positive element $y$  satisfying the relations \[ w_{1}w_{1}^{*}+...+ w_{m}w_{m}^{*}+ y^{2}=1, \quad w_{i}w_{j}=\rho_{ji}  w_{j}w_{i}, \quad [y, w_{i}]=0. \]
\end{definition}

\begin{prop}
\label{S4thetapullback}
The C*-algebra $C(S^{2m}_{\theta})$ is isomorphic to the pullback \[ \xymatrix{P \ar[d]_{i_2} \ar[r]^{i_1} & C(D^{2m}_{\theta}) \ar[d]^{j_1} \\ C(D^{2m}_{\theta}) \ar[r]_{j_2} & C(S^{2m-1}_{\theta}),}\] where the$i_{k}$'s and $j_{k}$'s are specified below.
\end{prop}
\begin{proof}
Let $A_{1}$ and $A_{2}$ be two copies of $C(D^{2m}_{\theta})$.
Consider the maps $i_{k}: C(S^{2m}_{\theta}) \rightarrow A_{1} \oplus A_{2}$ given by \[ i_{1}(z_{i})=(w_{i}, 0) \quad i_{1}(x)=(y, 0) \quad i_{2}(z_{i})=(0, w_{i}) \quad i_{2}(x)=(0, -y).\]

The maps $i_{k}$ are well-defined by the universal property of $C(S^{2m}_{\theta})$. Now let $i:C(S^{2m}_{\theta}) \rightarrow A_{1} \oplus A_{2}$ be the map $i=i_{1}+i_{2}$.

We wish to identify the image of the map $i$ in $A_{1} \oplus A_{2}$.  We let $j_{k}: A_{k} \rightarrow C(S^{2m-1}_{\theta})$ be the quotient map defined by $j_{k}(w_{i})=W_{i}$ and $j_{k}(y)=0$, where the $W_{i}$ generate $C(S^{2m-1}_{\theta})$, and define the pullback $P=\{(a_{1}, a_{2}) \in A_{1} \oplus A_{2} : j_{1}(a_{1})=j_{2}(a_{2})\}$. Now the image of the map $i$ is generated by the elements $(w_{i}, w_{i})$ and $(y, -y)$, which are in $P$. So $i(C(S^{2m}_{\theta})) \subset P$. Thus to see that $i(C(S^{2m}_{\theta}))$ is isomorphic to $P$, it suffices to see that the elements $(w_{i}, w_{i})$ and $(y, -y)$ generate $P$. Since the ideal $A_{k}y$ is contained in ${\mathrm ker}(j_k)$, and since the images of the generators of $A_{k}$ in the quotient $A_{k}/A_{k}y$ satisfy the relations in $C(S^{2m-1}_{\theta}) \cong A_{k}/{\mathrm ker}(j_k)$, we get that ${\mathrm ker}(j_k)=A_{k}y$, by universality of $C(S^{2m-1}_{\theta})$. From this and the fact that $A_{1}$ and $A_{2}$ are copies of each other, it is clear that $(a_{1}, a_{2}) \in P$ if and only if we can write $a_{1}= a+y_{1}$ and $a_{2}=a+y_{2}$, for some $a \in A_{1}=A_{2}$, some $y_{1} \in A_{1}y$ and some $y_{2} \in A_{2}y$. But the $*$- subalgebra of elements of $A_{1}\oplus A_{2}$ that can be written in such form is easily seen to be generated by the elements $(w_{i}, w_{i})$ and $(y, -y)$.

We wish to show that $i$ is injective. To this end we use the continuous functional calculus to write $x=x_{+}-x_{-}$, where $x_{+}x_{-}=x_{-}x_{+} = 0$, and $x_{+}$ and $x_{-}$ are positive. Since $(y, 0)=i_{1}(x_{+})-i_{1}(x_{-})$,  we see that the ideal $C(S^{2m}_{\theta})x_{-}$ is contained in $I_{1}={\mathrm ker}(i_1)$. So the map $C(S^{2m}_{\theta})/C(S^{2m}_{\theta})x_{-} \rightarrow C(S^{2m}_{\theta})/ I_{1} \cong A_{1}$ is surjective. On the other hand, the algebra $C(S^{2m}_{\theta})/C(S^{2m}_{\theta})x_{-}$ is generated by the images of the generators of $C(S^{2m}_{\theta})$ under the quotient map $C(S^{2m}_{\theta}) \rightarrow C(S^{2m}_{\theta})/C(S^{2m}_{\theta})x_{-}$. The images of the generators satisfy the relations that the $w_{i}$ and $y$ satisfy in $A_{1}$. So by universality of $A_{1}$ we get a surjection $A_{1} \rightarrow  C(S^{2m}_{\theta})/C(S^{2m}_{\theta})x_{-}$. Thus $I_{1}=C(S^{2m}_{\theta})x_{-}$. By a similar argument, we also have that $I_{2}={\mathrm ker}(i_2)=C(S^{2m}_{\theta})x_{+}$. Now since ${\mathrm image}(i_1)\cap {\mathrm image}(i_2)=(0,0)$, we have that ${\mathrm ker}(i_{1}+i_{2}) \subset I_{1} \cap I_{2}$. Since $I_{2}$ has an approximate identity, the ideal $I_{1}I_{2}$ is all of $I_{1} \cap I_{2}$. But since $x_{+}$ and $x_{-}$ commute with all elements of $C(S^{2m}_{\theta})$, and since $x_{+}x_{-}=0$, the ideal $I_{1}I_{2}$ is just $\{0\}$.

\end{proof}

Note that $C(S^{2m-1}_{\theta})$ is not the pullback of two $2m-1$-dimensional $\theta$-deformed balls, but does respect decompositions of $S^{2m-1}$ that are equivariant for the canonical action of $T^{m}$ on $S^{2m-1}$.

\begin{theorem}
\label{evenfield}
Let $S(i)=\{(t_{1},..., t_{m+1}) \in S^{m}_{+} : t_{i}=0\}$. Let $A_{\theta}(i)$ be the C*-subalgebra of $A_{\theta}$ generated by the set $\{U_{1},...U_{m}\}$ omitting $U_{i}$. Let $\mathbb{A}$ be the C*-subalgebra of functions in $C(S^{m}_{+}, A_{\theta})$ which map $S(i)$ into $A_{\theta}(i)$ for all $ i\leq m$. Then $C(D^{2m}_{\theta})\cong \mathbb{A}$.
\end{theorem}
\begin{proof}
Our argument is a fairly straightforward adaptation of the proof of theorem 2.5 of \cite{natsume}. For much of the proof, we follow \cite{natsume} nearly verbatim.

 It is clear that the functions \[(t_{1},..., t_{m+1}) \mapsto t_{i}U_{i} \quad i=1,..., m \] represent the $w_{i}$, and that the function \[ (t_{1},..., t_{m+1}) \mapsto \sqrt{1-t_{1}^{2}-...-t_{m}^{2}} \cdot 1\] represents $y$. 

That $\mathbb{A}$ is generated as a C*-algebra by these functions is proved along the lines of proposition 2.4 of \cite{natsume}.

Thus it suffices to show that if a C*-algebra $\mathbb{B}$ is generated by $m$ normal elements $B_{1},..., B_{m}$, and a positive element $B_{m+1}$, satisfying the relations  \[ B_{1}B_{1}^{*}+...+ B_{m}B_{k}^{*}+ B_{m+1}^{2}=1, \quad B_{i}B_{j}=\rho_{ji}  B_{j}B_{i}, \text{ for all }i, j \leq m,\] and \[ [B_{m+1}, B_{i}]=0,\text{ for all }i \leq m+1,\] that there is then a surjection $\pi$ from $\mathbb{A}$ onto $\mathbb{B}$. 

As in \cite{natsume}, we first prove the existence of the surjection $\pi$ under the assumption that all of the $B_{i}$ have finite spectra, and then use an approximation argument to handle the general case.

Under the assumption that the $B_{i}$ have finite spectra we have that the partial isometries in the polar decompositions of the $B_{i}$ are in $\mathbb{B}$. This will be used in the proof.

Since the $B_{i}$ have finite spectra, we may write $|B_{i}|=\Sigma \mu_{i k}E_{i,k}$ for each $i \leq m+1$, where $\{E_{i, k}\}_{k=1,...m_{i}}$ is a pairwise orthogonal set of spectral projections of $|B_{i}|$. We include the projection null $B_{i}$ in $\{E_{i,k}\}_{k}$ so that $\Sigma_{k}E_{i, k}=I$. 

As $B_{i}$ is normal, $B_{i}= \Sigma_{k}E_{ik}B_{i}E_{ik}$, and $E_{ik}B_{i}E_{ik}=B_{i}E_{ik}= \mu_{ik}V_{ik}E_{ik}$, where if $E_{ik}=$null$|B_{i}|$, we choose $V_{ik}=0$, and otherwise $V_{ik}^{*}V_{ik}=E_{ik}=V_{ik}V_{ik}^{*}$.

 Set $V_{i}= \Sigma_{k}V_{ik}$, so that $V_{i}E_{ik}=V_{ik}E_{ik}=V_{ik}$, and $B_{i}=V_{i}|B_{i}|$ is the polar decomposition of $B_{i}$. As in \cite{natsume}, we have that $V_{j}V_{i}= \rho_{ij} V_{i}V_{j}$ for all $i, j \leq m$, and we also have that $V_{j}V_{m+1}=V_{m+1}V_{j}$ for all $j \leq m+1$.

 Now consider the set $\{E_{1k_{1}}...E_{(m+1)k_{m+1}}\}_{k_{i}=1}^{m_{i}},$ $i=1,..., m+1$. This set of projections is pairwise orthogonal, and after throwing away any projections which are zero,  we renumber the set as $P_{1},...P_{d}$. Each $P_{j}$ is in $C^{*}(B_{1},..., B_{m+1})$ and $\Sigma_{j} P_{j}=I$. 

As in \cite{natsume}, each $B_{i}$ commutes with every $P_{j}$ and $B_{i}=\Sigma_{j}P_{j}B_{i}P_{j}$, so $\mathbb{B}$ is a finite orthogonal direct sum of the algebras $P_{j}B$.

 Each $B_{i}P_{j}=P_{i}B_{i}P_{j}$ has the form $B_{i}P_{j}= \lambda_{ij}V_{ij}P_{j}$, where $V_{ij}$ is defined as $P_{j}V_{i}P_{j}$ (notice that as in \cite{natsume}, this notation is not consistent with our earler use of $V_{ik}$) and $\lambda_{ij} \geq 0$. This $\lambda_{ij}$ is one of $\mu_{i1},...\mu_{im_{i}}$, and by our choice, if $\lambda_{ij}=0$, then $V_{ij}=0$. If $V_{ij} \neq 0$, then $V_{ij}$ is unitary on $P_{j}\mathbb{B}$.

We have that for each $i$, \[ |B_{i}|= \Sigma_{j} \lambda_{ij}P_{j}, \quad |B_{i}|^{2} = \Sigma_{j} \lambda_{ij}^{2}P_{j}, \] so that for any fixed $j$, we have \[ P_{j}= (\Sigma_{i}|B_{i}|^{2})P_{j}= \Sigma_{i}(\Sigma_{k} \lambda_{ik}^{2})P_{j} = \Sigma_{i} \lambda_{ij}^{2}P_{j}. \] Thus $\Sigma_{i} \lambda_{ij}^{2}=1$, for each fixed $j$, and so the point $(\lambda_{1j},...\lambda_{(m+1)j})$ is in $S^{m}_{+}$.

We fix an arbitrary $j$ for what follows, and restrict to the summand $P_{j}\mathbb{B}$ (on which each $B_{i}P_{j}= \lambda_{ij} V_{ij}$, and $V_{ij}$ is either unitary on $P_{j}\mathbb{B}$ or $V_{ij}=0$). Then $B_{i}B_{k}P_{j}= \lambda_{ij}\lambda_{kj}V_{ij}V_{kj}$, but also $B_{i}B_{k}P_{j}= \rho_{ki}B_{k}B_{i}P_{j} = \rho_{ki} \lambda_{ij} \lambda_{kj} V_{kj}V_{ij}$, so we have the commutation relations $V_{ij}V_{kj}= \rho_{ki}V_{kj}V_{ij}$, for $i, k = 1,..., m$, and $V_{ij}V_{(m+1), j}= V_{(m+1), j}V_{ij}$ for $i=1,..., m+1$.

We now consider the different possibilities for the values of the $\lambda_{ij}, i=1,..., m+1$, showing how to define a surjective $*$-homomorphism $\pi_{j}:\mathbb{A} \rightarrow P_{j}\mathbb{B}$ in each case. So far we have been following the argument of \cite{natsume} closely. It is in defining these surjections that we need to modify the argument of \cite{natsume} somewhat nontrivially.

First suppose that $\lambda_{(m+1)j}$ is not zero, and exactly $k< m$ of the $\lambda_{1j},..., \lambda_{mj}$ are zero, namely $\lambda_{i_{1}j},... \lambda_{i_{k}j}$. Then $V_{i_{1}j}=...V_{i_{k}j}=0$, and the rest of the $V_{ij}$ are unitary. Let $A^{m-k}_{{\theta}^{\prime}}$ be the C*-subalgebra  of $A^{m}_{\theta}$ generated by $U_{i} \in A^{m}_{\theta}$ excluding the $k$ unitaries $U_{i_{1}},..., U_{i_{k}}$, where $\theta^{\prime}$ is the $m-k \times m-k$ matrix formed from $\theta$ by deleting the rows and columns corresponding to $i_{1},..., i_{k}$. The C*-algebra $A^{m-k}_{{\theta}^{\prime}}$ is itself universal, and so there is a  $*$-homomorphism $A^{m}_{\theta} \rightarrow P_{j}\mathbb{B}$ mapping $U_{i}$ to $V_{ij}$ for $ i \neq i_{1},...,i_{k}$. By the functional calculus, there is also a $*$-homomorphism $C([0,1]) \rightarrow P_{j} \mathbb{B}$ mapping the generator $f(t)=t$ of $C([0, 1])$ to $B_{m+1}P_{j}$. Since  $C([0,1])$ is  nuclear, there is a unique C*-cross norm on the algebraic tensor product $C([0,1]) \otimes_{alg} A^{m-k}_{{\theta}^{\prime}}$. Completing in this norm we obtain a C*-algebra $C([0,1]) \otimes A^{m-k}_{{\theta}^{\prime}}$. Since $B_{m+1}P_{j}$ commutes with all of the $V_{ij}$, and being that the tensor product $C([0,1]) \otimes A^{m-k}_{{\theta}^{\prime}}$ is maximal, we get a surjective $*$-homomorphism $\phi_{j}: C([0,1]) \otimes A^{m-k}_{{\theta}^{\prime}} \rightarrow P_{j}\mathbb{B}$ by the universal property of the tensor product. We can then define $\pi_{j}$ as follows: 

Consider the set 
\begin{align*} 
J=\{(s_{1},...s_{m+1}) \in S^{m}_{+} : \quad s_{i} \geq \lambda_{ij} \quad i \neq i_{1},...,i_{k}, \quad i \leq m, \\
 \text{and } \frac{s_{p}}{s_{q}}=\frac{\lambda_{pj}}{\lambda_{qj}} \quad  p,q \neq i_{1},...,i_{k}, \quad p,q \leq m, \\
 \text{and } s_{i_{1}}=...=s_{i_{k}}= 0 \}.
\end{align*}

Clearly $J$ is homeomorphic to $[0,1]$ via an obvious map that takes the point $(\lambda_{1j},...\lambda_{(m+1)j})$ to $1$ and takes the unique point in $J$ of the form $(s_{1},...s_{m}, 0)$ to $0$. So given $f \in \mathbb{A}$, we may restrict $f$ to the set $J$, and then identify this restriction with a map $C([0,1]) \otimes A^{m-k}_{{\theta}^{\prime}}$ by using the homeomorphism $J \approx [0,1]$. We then compose with the map $\phi_{j}$ to obtain a surjective $*$-homomorphism $\pi_{j}:\mathbb{A}\rightarrow P_{j}B$.

The next case we consider is when $\lambda_{(m+1)j}$ is zero, and as before exactly  $k< m$ of the $\lambda_{1j},..., \lambda_{mj}$ are zero, specifically $\lambda_{i_{1}j},... \lambda_{i_{k}j}$. Then we form $A^{m-k}_{{\theta}^{\prime}}$ as in the previous case, and now since we have in addition that $P_{j}B_{(m+1)j}$ is zero,  the $*$-homomorphism $A^{m-k}_{{\theta}^{\prime}} \rightarrow P_{j}\mathbb{B}$ is surjective. By composing this map with evaluation at the point $(\lambda_{1j},...\lambda_{(m+1)j})$,  we obtain a surjective $*$-homomorphism $\pi_{j}:\mathbb{A} \rightarrow P_{j}\mathbb{B}$.

The remaining case is when $\lambda_{ij}$ is zero for all $i \leq m$, so that $\lambda_{(m+1)j}$ equals one. In this case, we have that $B_{i}P_{j}$ is zero for all $i \leq m$, and $B_{(m+1)j}$ is the identity operator on $P_{j}\mathbb{B}$. Thus we have a surjective $*$-homomorphism $\mathbb{C} \rightarrow P_{j}\mathbb{B}$, and so again composing this map with evaluation at the point $(\lambda_{1j},...\lambda_{(m+1)j})$,  we obtain a surjective $*$-homomorphism $\pi_{j}:\mathbb{A} \rightarrow P_{j}\mathbb{B}$.

Define the surjection $\pi: \mathbb{A} \rightarrow \mathbb{B}$ by $\pi(a)= \Sigma_{j} \pi_{j}(a)$, for each $a \in A$. Since the ranges of the $\pi_{j}$ are pairwise orthogonal, the surjection $\pi$ is a $*$-homomorphism.

One finishes the proof by approximating to the case where the $B_{i}$ have infinite spectra as in \cite{natsume}.

\end{proof}

Theorem \ref{evenfield} generalizes the decomposition of $C(D^{2m})$ as the fibered space over $S^{m}_{+}$ that comes from the action of $T^{m}$ on $S^{2m}$ that rotates the first $m$ coordinates of a vector in $S^{2m}$, but leaves the $m+1$-st coordinate fixed.

Let \[S^{m}_{+, m+1} = \{(t_{1},..., t_{m+1}) \in S^{m} : t_{k} \geq 0  \text{ for all } k \leq m \}.\] The space $S^{m}_{+, m+1}$ is clearly homeomorphic to the union of two $m$-simplices glued along one face.

From Theorem \ref{evenfield} and Proposition \ref{S4thetapullback} we immediately obtain:

\begin{cor}
\label{evenspherefield}
 Let $S(i)=\{(t_{1},..., t_{m+1}) \in S^{m}_{+, m+1} : t_{i}=0\}$. Let $A_{\theta}(i)$ be the C*-subalgebra of $A_{\theta}$ generated by the set $\{U_{1},...U_{m}\}$ omitting $U_{i}$. Let $\mathfrak{A}$ be the C*-subalgebra of functions in $C(S^{m}_{+}, A_{\theta})$ which map $S(i)$ into $A_{\theta}(i)$ for all $ i\leq m$. Then $C(S^{2m}_{\theta})\cong \mathfrak{A}$.
\end{cor}

Let $A^{\sim}$ denote the unitisation of a C*-algebra $A$.

\begin{cor}
\label{cone} Let $m \geq 2$. Then
\begin{equation*}
\begin{split}
C(D^{2m}_{\theta}) & \cong \{ f \in C([0,1], C(S^{2m-1}_{\theta})) : f(1) \in \mathbb{C}\cdot 1 \}\\
&  \cong (C_{0}([0,1)) \otimes C(S^{2m-1}_{\theta})^{\sim}
\end{split}
\end{equation*}
\end{cor}
\begin{proof}
Let $\mathbb{A} \cong C(D^{2m}_{\theta})$ be as in Theorem \ref{evenfield}. For fixed $t \in [0, 1]$, we may consider the C*-algebra of all restrictions of functions in $\mathbb{A}$ to the set $\{(t_{1},..., t_{m+1}) \in S^{m+1}_{+} : t_{1}^{2}+...t_{m}^{2}= t^{2} \}$. By Theorems \ref{oddfield} and \ref{evenfield}, this C*-algebra is isomorphic to $C(S^{2m-1}_{\theta})$ if $t \neq 0$, but is isomorphic to $\mathbb{C}$ if $t = 0$. Thus by varying $t$ and reversing the parametrization of the interval $[0,1]$, we obtain that $\mathbb{A} \cong \{ f \in C([0,1], C(S^{2m-1}_{\theta})) : f(1) \in \mathbb{C} \cdot 1 \}$. The second isomorphism is clear from the identity $f= f-f(1)+f(1)$.
\end{proof}

\begin{cor}
\label{deformationquant}
Let $l=2m$ or $2m-1$, and let $\alpha$ be the canonical $T^{m}$ action on $C(S^{l})$ viewed as an action of $\mathbb{R}^{m}$. Then $C(S^{l}_{\theta})$ is a deformation quantization (in the sense of \cite{deformation}) of $C(S^{l})$ by the $\mathbb{R}^{m}$-action $\alpha$.
\end{cor}
\begin{proof}
The case $l=2m-1$ is proposition 2.9 of \cite{natsume}. In view of Corollary \ref{evenspherefield}, one sees that the proof of proposition 2.9 of \cite{natsume} also applies in the case $l=2m$.
\end{proof}

The $K$-groups of the $C(S^{l}_{\theta})$ are known:

\begin{theorem}
\label{K-groups}
 \[ K_{0}(C(S^{2m-1}_{\theta}))\cong \mathbb{Z}, \quad K_{1}(C(S^{2m-1}_{\theta}))\cong \mathbb{Z}.\]
\[ K_{0}(C(S^{2m}_{\theta}))\cong \mathbb{Z} \times \mathbb{Z}, \quad K_{1}(C(S^{2m}_{\theta}))\cong 0.\]
\end{theorem}

Matsumoto \cite{matsumoto} proved Theorem \ref{K-groups} for the case $l=3$ by using the Mayer-Vietoris sequence for $K$-theory. The general odd case is shown in \cite{natsume} by using a Schochet spectral sequence \cite{schochet}. As noted in \cite{natsume}, it follows from this together with Theorem \ref{oddfield}, that if $\theta$  is irrational, then a generator for $K_{1}(C(S^{3}_{\theta}))$ is given explicitly by $X=\exp(2\pi i t)p + 1-p$, where $p$ is a Rieffel projection of trace $\theta \mod 1$, and $S^{1}_{+}$ is identified with the unit interval via the composition $S^{1}_{+} \approx \Delta^{1} \approx [0,1]$. We observe in Theorem \ref{naturalmap} that in fact $X$ generates $\pi_{0}(GL_{n}(C(S^{3}_{\theta})))$ for each $n \geq 1$. This observation will be fundamental to our classification and construction of all finitely-generated projective $C(S^{4}_{\theta})$-modules up to isomorphism.

The general even case of Theorem \ref{K-groups} follows easily from \cite{natsume} by applying the $K$-theory Mayer-Vietoris sequence to the pullback given in Proposition \ref{S4thetapullback}, and noting that $C(D^{4}_{\theta})$ deformation retracts onto $\mathbb{C}$ (see Proposition \ref{D4trivial}), or by citing Corollary \ref{deformationquant} together with the main theorem of \cite{kgroups}.
\newline

Looking for examples of noncommutative differentiable manifolds in the sense of \cite{gravity}, Connes and Landi \cite{isospectral} found the spheres $C(S^{l}_{\theta})$ as solutions to certain homological equations, and also showed that they could be seen as coming from a type of deformation procedure they termed ``$\theta$-deformation". It was soon realized \cite{sitarz, varilly} that $\theta$-deformation in general is a special case of Rieffel's deformation quantization by actions of $\mathbb{R}^{m}$. 

Connes and Dubois-Violette found all solutions of (the purely algebraic form of) Connes and Landi's equations for $l=3$ as a three-parameter family of algebras in \cite{connesdv}. For critical values in the moduli space of solutions, one obtains algebras whose universal C*-completions are isomorphic to the $C(S^{3}_{\theta})$. In the generic case, one recovers the Sklyanin algebras \cite{sklyanin} of noncommutative algebraic geometry. Connes and Dubois-Violette went on to investigate the structure of the moduli space of the three-parameter family of algebras in \cite{moduli}, exhibiting very interesting further connections to algebraic geometry and Morse theory. Building on the preliminary consideration of the Yang-Mills theory on $C(S^{4}_{\theta})$ initiated in \cite{isospectral}, Brain, Landi and van Suijelkom (\cite{brainvs}, \cite{brain}) relate moduli of instantons on $C(S^{4}_{\theta})$ to algebraic bundles on $\theta$-deformed twistor space using a modified ADHM construction. One wonders whether, in addition, there may be relationships to string theory or to the moduli of classical algebraic varieties as for the modified ADHM equations studied in \cite{nekrasov}, \cite{kapustin}, \cite{nakajima}.

\section{Clutching Construction}
Throughout this section we work in the category of unital rings.
\newline
Our classification of the right finitely-generated projective $C(S^{3}_{\theta})$ and $C(S^{4}_{\theta})$-modules relies on a generalization of the familiar clutching construction (theorem 3.1 of \cite{karoubi}) of vector bundle theory to the setting of unital rings. We follow Milnor \cite{milnor}, and the reader familiar with the results there may skip this section, referring back to it as needed.

Given $ M$ a right $ R$-module, a unital ring homomorphism $ k:R\rightarrow S$ induces a right $S$-module $k_{\#}M:=M\otimes_{R}S$  and a canonical $R$-linear map $k_{*}:M\rightarrow k_{\#}M$ defined by $k_{*}(m):=m\otimes_{R}1$.
This process is functorial, so that if $g: M \rightarrow M^{\prime}$ is an $R$-module map, there is then an $S$-module map $k_{*}(g): k_{\#}M \rightarrow k_{\#}M^{\prime}$ given by $k_{*}(g)(m \otimes_{R}1)=g(m) \otimes_{R} 1$.
 
Consider the pullback diagram \[ \xymatrix{R \ar[d]_{i_2} \ar[r]^{i_1} & R_{1} \ar[d]^{j_1} \\ R_{2} \ar[r]_{j_2} & S.}\]
Now suppose $P_{1}$ is a right $R_{1}$-module, $P_{2}$ is a right $R_{2}$-module, and that $h:j_{1\#}P_{1} \rightarrow j_{2\#}P_{2}$ is an isomorphism. 
\begin{definition}
\label{module}
 The right $R$-module $M(P_{1}, P_{2},h)$ is the additive group \[\{(p_{1},p_{2}) \in P_{1} \times P_{2} : hj_{1*}(p_{1})=j_{2*}(p_{2})\}\] together with the $R$-module structure given by \[(p_{1}, p_{2}) \cdot r = (p_{1} \cdot i_{1}(r), p_{2} \cdot i_{2}(r)).\]
\end{definition}
Assume also that at least one of the $j_{k}$ is surjective. Then we have the following three theorems:
\begin{theorem}
\label{miln1}
(Theorem 2.1 of \cite{milnor}). If $P_{1}$ and $P_{2}$ are projective over $R_{1}$ and $R_{2}$ respectively, then $M(P_{1}, P_{2}, h)$ is projective over $R$. If also $ P_{1}$ and $P_{2}$ are finitely-generated, then $ M(P_{1}, P_{2}, h)$ is finitely-generated.
\end{theorem}

\begin{theorem}
\label{miln2}
(Theorem 2.2 of \cite{milnor}). Any projective $R$-module $M$ is isomorphic to $M(P_{1}, P_{2}, h)$ for some $P_{1}, P_{2}$ and $ h$.
\end{theorem}

\begin{theorem}
\label{miln3}
(Theorem 2.3 of \cite{milnor}). The modules $P_{1}$ and $P_{2}$ are canonically isomorphic to $i_{1\#}M$ and $i_{2\#}M$ respectively.
\end{theorem}

We will use the following corollary to Theorem \ref{miln3} in classifying all finitely-generated projective $C(S^{4}_{\theta})$-modules up to isomorphism.

\begin{cor}
\label{milncor}
Two $R$-modules $M=M(P_{1}, P_{2}, h)$ and $N=M(Q_{1}, Q_{2}, g)$ are isomorphic if and only if there exist isomorphisms $g_{k}: Q_{k} \rightarrow P_{k}$ such that $h=j_{2*}(g_{2})g j_{1*}(g_{1}^{-1})$.
\end{cor}
\begin{proof}
If $h=j_{2*}(g_{2})g j_{1*}(g_{1}^{-1})$ where $g_{k}: Q_{k} \cong P_{k}$, then $M(Q_{1}, Q_{2}, g) \cong M(P_{1}, P_{2}, h)$ via the map $(q_{1}, q_{2}) \mapsto (g_{1}q_{1}, g_{2}q_{2})$.

To prove the converse, let $m=(p_{1}, p_{2})$. By Theorem \ref{miln3}, the natural map $i_{k\#}M \rightarrow P_{k}$ given by $m \otimes_{R} 1 \mapsto p_{k}$ is an isomorphism.
Thus we have a commutative diagram
\[ \xymatrix{ P_{1} \ar[d] & i_{1 \#}M  \ar[l]^{\cong} \ar[d] \\
j_{1\#}P_{1} \ar[d]_{h} & j_{1\#} i_{1 \#}M  \ar[l]^{\cong} \ar[d]^{f_{M}} \\
j_{2\#}P_{2} &  j_{2\#}i_{2\#}M \ar[l]^{\cong} \\
 P_{2} \ar[u] &  i_{2 \#}M ,\ar[l]^{\cong} \ar[u] } \]
where $f_{M}$ is the canonical isomorphism \[f_{M}(m\otimes_{R} 1 \otimes_{R_{1}} 1) = m \otimes_{R} 1 \otimes_{R_{2}} 1, \quad m \in M.\]
In the same way, we have
\[ \xymatrix{ Q_{1} \ar[d] & i_{1 \#}N  \ar[l]^{\cong} \ar[d] \\
j_{1\#}Q_{1} \ar[d]_{g} & j_{1\#} i_{1 \#}N  \ar[l]^{\cong} \ar[d]^{f_{N}} \\
j_{2\#}Q_{2} &  j_{2\#}i_{2\#}N \ar[l]^{\cong} \\
 Q_{2} \ar[u] &  i_{2 \#}N .\ar[l]^{\cong} \ar[u] } \]
Now if $N \cong M$, this induces isomorphisms $i_{k\#}(N) \rightarrow i_{k\#}(M)$, and so we have a commutative diagram
\[ \xymatrix{ i_{1 \#}M \ar[d] & i_{1 \#}N  \ar[l]^{\cong} \ar[d] \\
j_{1\#}i_{1\#}M \ar[d]_{f_{M}} & j_{1\#} i_{1 \#}N  \ar[l]^{\cong} \ar[d]^{f_{N}} \\
j_{2\#}i_{2\#}M &  j_{2\#}i_{2\#}N \ar[l]^{\cong} \\
 i_{2 \#}M \ar[u] &  i_{2 \#}N .\ar[l]^{\cong} \ar[u] } \]
Combining the three diagrams and composing along the second and third rows gives the maps $g_{1}$ and $g_{2}$.
\end{proof}

An immediate consequence of Theorems \ref{miln2}, \ref{miln3}, and Corollary \ref{milncor} is that any projective $R$-module $ M$ is isomorphic to $M(i_{1\#}M, i_{2\#}M, f_{M})$, where $f_{M}$ is the canonical isomorphism \[f_{M}(m\otimes_{R} 1 \otimes_{R_{1}} 1) = m \otimes_{R} 1 \otimes_{R_{2}} 1, \quad m \in M.\] This is the algebraic analog of the fact that a vector bundle $E$ over a compact space $X=Y \cup Z$ ($Y$ and $Z$ closed) can be viewed as $E |Y$ and $E|Z$ glued together over $Y \cap Z$ by the identity map.

The relationship $h=j_{2*}(g_{2})g j_{1*}(g_{1}^{-1})$ is an equivalence relation, so by Corollary \ref{milncor} there is a bijection between the set of equivalence classes of this equivalence relation and the set of isomorphism classes of finitely-generated projective $R$-modules. In general, it is very difficult to determine these equivalence classes. But if the $R_{k}$'s and $S$ satisfy the invariance of dimension property, then we may immediately draw some conclusions that simplify the situation somewhat.

\begin{definition}
A unital ring $R$ has the \emph{invariance of dimension} property if for each free $R$-module $F$, every basis for $F$ has the same cardinality.
\end{definition}
Invariance of dimension is equivalent to the condition that for any $n$ and $m$, the $R$-modules $R^{n}$ and $R^{m}$ are isomorphic if and only if $n=m$. It is easy to see that if $S$ has the invariance of dimension property, and $j:R \rightarrow S$ is surjective, then $R$ also has the invariance of dimension property. Any unital commutative ring has the invariance of dimension property.

\begin{prop}
\label{invtech1}
Suppose that $R_{1}$ and $R_{2}$ satisfy the invariance of dimension property. Then it is impossible for $M=M(R_{1}^{n}, R_{2}^{n}, g)$ to be isomorphic to $N= M(R_{1}^{m}, R_{2}^{m}, h)$, unless $n=m$.
\begin{proof}
By Theorem \ref{miln3} we have that $i_{k\#}M \cong R_{k}^{n}$ and $i_{k\#}N \cong R_{k}^{m}$. So $M$ cannot be isomorphic to $N$ unless $n=m$, since the rings $R_{k}$ satisfy the invariance of dimension property.
\end{proof}
\end{prop}

\begin{prop}
\label{invtech2}
Suppose that $R_{1}$ and $R_{2}$ satisfy the invariance of dimension property. Then a $R$-module can be of the form $M(R_{1}^{n}, R_{2}^{m}, h)$ only if $n =m$.
\begin{proof}
We have that $j_{1\#}R_{1}^{n} \cong S^{n}$, and also that $j_{2\#}R_{1}^{m} \cong S^{m}$. So if $h$ is an isomorphism, then $S^{n} \cong S^{m}$, which is possible only if $n=m$, since $S$ satisfies the invariance of dimension property.
\end{proof}
\end{prop}

We note that if $\psi_{k}:P_{k} \cong R_{k}^{n}$, then, by Corollary \ref{milncor}, we may replace $M(P_{1}, P_{2}, h)$ with $M(R_{1}^{n}, R_{2}^{n}, j_{2*}(\psi_{2})hj_{1*}(\psi_{1}^{-1}))$. Additionally, we may identify $Aut_{R_{k}}(R_{k}^{n})$ with $GL_{n}(R_{k})$  by choosing a basis for $R_{k}^{n}$. Under such an identification, an automorphism of $R^{n}_{k}$ acts on elements of $R^{n}_{k}$ (viewed as column vectors) by left multiplication (since we are dealing with right $R_{k}$-modules) by the corresponding matrix in $GL_{n}(R_{k})$. Throughout this work, we always choose the standard ordered basis $\{e_{k,i}\}$ for $R_{k}^{n}$. Now $j_{k\#}R_{k}^{n} \cong S^{n}$, via the map that takes $e_{k,i}\otimes_{R_{k}} 1$ to the $i$-th element of the standard ordered basis for $S^{n}$. Thus, if $h$ is an isomorphism $j_{1\#}R_{1}^{n} \cong j_{2\#}R_{2}^{n}$, then we may identify $h$ with a matrix in $GL_{n}(S)$ that acts on elements of  $S^{n}$ viewed as column vectors by left multiplication. If $h$ is the isomorphism given by $h(e_{1,i} \otimes_{R_{1}} 1)=e_{2, i} \otimes_ {R_{2}} 1$, then $M(R_{1}^{n}, R_{2}^{n}, h)$ is isomorphic to $R^{n}$, and $h$ identifies with the identity matrix $1_{n} \in GL_{n}(S)$. Therefore we will refer to the map $h(e_{1,i} \otimes_{R_{1}} 1)=e_{2, i} \otimes_ {R_{2}} 1$ as $id$.

As short-hand for $M(R_{1}^{n}, R_{2}^{n}, h)$, we will occasionally write expressions of the form $M(R_{1}^{n}, R_{2}^{n}, \alpha)$, where $\alpha \in GL_{n}(S)$ is the matrix corresponding to $h$ under our identification of the isomorphism group $Iso_{S}( j_{1\#}R^{n}_{1}, j_{2\#}R^{n}_{2})$ with the matrix group $GL_{n}(S)$.

\section{Finitely-Generated Projective Modules over $C(S^{3}_{\theta})$}

Throughout this section $\theta$ is an arbitrary real number unless otherwise specified. We prove that all right finitely-generated projective $C(S^{3}_{\theta})$-modules are free, and that $V(C(S_{\theta}^{3})) \cong \mathbb{N}$.

Our strategy is to use the clutching construction of Section 2. We think of $C(S^{3}_{\theta})$ as consisting of two noncommutative solid tori ``hemispheres"  with a noncommutative torus as the ``equator" that we clutch over.

We first need a series of small propositions. Note that we now work in the category of unital C*-algebras (of course this restriction could be considerably weakened for our results), rather than in the purely algebraic category of unital rings.

Let \[\mathcal{A}=  \{f \in C([0,1], C(T^{2}_{\theta})) : f(0) \in C^{*}(U_{2}), f(1) \in C^{*}(U_{1})\},\] \[\mathcal{A}_{1}= \{f \in C([1/2,1], C(T^{2}_{\theta})) : f(1) \in C^{*}(U_{1})\},\] \[\mathcal{A}_{2}= \{f \in C([0,1/2], C(T^{2}_{\theta})) : f(0) \in C^{*}(U_{2}) \},\] and \[ \mathcal{B}= \{f \in C(\{1/2\}, C(T^{2}_{\theta})) \}.\]

\begin{prop}
\label{3spherepullback}
The C*-algebra $C(S^{3}_{\theta})$ is isomorphic to the pullback  \[ \xymatrix{ \mathcal{A}\ar[d]_{i_2} \ar[r]^{i_1} &  \mathcal{A}_{1}  \ar[d]^{j_1} \\  \ \mathcal{A}_{2} \ar[r]_{j_2} & \mathcal{ B},}\] 
where the $i_{k}$'s and $j_{k}$'s are given by restriction.
\end{prop}
\begin{proof}
That $\mathcal{A}$ is given by the pullback diagram is clear. That $C(S^{3}_{\theta})$ is isomorphic to $ \mathcal{A}$ follows from Theorem \ref{oddfield}, as explained immediately after the statement of that theorem.
\end{proof}

\begin{prop}
\label{invdimension}
 The C*-algebras $\mathcal{B}$, $\mathcal{A}_{1}$, $\mathcal{A}_{2}$, and $\mathcal{A}$ all have the invariance of dimension property.
\end{prop}
\begin{proof}
Choose a tracial state $Tr$ on $\mathcal{B} \cong C(T^{2}_{\theta})$ (note that if $\theta$ is irrational then $C(T^{2}_{\theta})$ has a unique tracial state, but there are many choices if $\theta$ is rational). Composing $Tr$ with the surjections  $j_{k}: \mathcal{A}_{k} \rightarrow \mathcal{B}$ and $j_{1}i_{1}: \mathcal{A} \rightarrow \mathcal{B}$ show that the $\mathcal{A}_{k}$'s and $\mathcal{A}$ also have (non-unique) tracial states. But any tracial state on a unital C*-algebra $A$ induces a semigroup homomorphism $V(A) \rightarrow \mathbb{R_{+}}$ that takes the class of the free $A$-module $A^{n}$ to $n$.
\end{proof}

\begin{prop} 
\label{filledmodstrivial}
All right finitely-generated projective $\mathcal{A}_{k}$-modules are free, for $k=1,2$.
\end{prop}
\begin{proof}
 We observe that $\mathcal{A}_{1}$ deformation retracts onto its commutative C*-subalgebra of constant $C^{*}(U_{1})$-valued functions. A retraction is given explicitly by $f_{t}(s) = f((1-t)s+ t)$. But the retracted C*-subalgebra is isomorphic to $C^{*}(U_{1}) \cong C(S^{1})$. The composition of this isomorphism with the retraction induces a semigroup homomorphism $V(\mathcal{A}_{1}) \rightarrow V(C(S^{1}))$. By the homotopy invariance of the functor $V$, this induced map is a isomorphism. Thus, since all finitely-generated projective modules over $C(S^{1})$ are free, the same is true for for all finitely generated projective modules over $\mathcal{A}_{1}$. Similar considerations apply to $\mathcal{A}_{2}$.
\end{proof}

Recall that a finitely-generated projective module $P$ over a unital C*-algebra $A$ may be given a norm so that the action of $A$ on $P$ is continuous in the corresponding topology, and so that $P$ is complete. This is because $P$ is the image of a free $A$-module $A^{n}$ under a projection, and $A^{n}$ inherits a norm naturally from the norm on $A$. The norm on $P$ does not depend on the choice of $n$ or on the choice of the projection, since any two such choices of projections are similar if we embed the projections into a large enough matrix algebra over $A$.  Of course the identification $Aut_{A}(A^{n}) \cong GL_{n}(A)$ of Section 2 now becomes one of topological groups.

\begin{prop}
\label{clutching}
 Suppose that \[ \xymatrix{A \ar[d]_{i_2} \ar[r]^{i_1} & A_{1}  \ar[d]^{j_1} \\ A_{2} \ar[r]_{j_2} & B}\] is a pullback in the category of unital  C*-algebras. Let $P_{1}$ be a projective $A_{1}$-module, and $P_{2}$  a projective $A_{2}$-module. Suppose that $h_{0}$ and $h_{1}$ are path-connected in $Iso_{B}(j_{1\#}P_{1}, j_{2\#}P_{2})$. Then $M(P_{1}, P_{2}, h_{0})$ and $M(P_{1}, P_{2}, h_{1})$ are isomorphic.
\end{prop}
\begin{proof}
 Form the induced pullback diagram  \[ \xymatrix{C(I, A) \ar[d]_{i_2} \ar[r]^{i_1} & C(I, A_{1})  \ar[d]^{j_1} \\ C(I,  A_{2}) \ar[r]_{j_2} & C(I, B).}\] We regard a path $\{h_{t}\}$ from $h_{0}$ to $h_{1}$ through $Iso_{B}(j_{1\#}P_{1}, j_{2\#}P_{2})$ as an element of $Iso_{C(I, B)}(j_{1\#}C(I,P_{1}), j_{2\#}C(I, P_{2}))$ and form the $C(I, A)$-module \[E=M(C(I, P_{1}),C(I,  P_{2}), \{h_{t}\}).\] Let $\varepsilon_{t}: C(I, A) \rightarrow A$ be evaluation at $t$ and apply the functor $V$ to obtain  a map $\varepsilon_{t*}: V(C(I,A))\rightarrow V(A)$. By homotopy invariance of the functor $V$, we have that $\varepsilon_{0*}$ and $\varepsilon_{1*}$ are equal. But simply by inspecting their definitions, we see that $\varepsilon_{0*}(E) \cong M(P_{1}, P_{2}, h_{0})$ and $ \varepsilon_{1*}(E) \cong M(P_{1}, P_{2}, h_{1})$.
\end{proof}

We identify isomorphisms $j_{1\#}C(T^{2}_{\theta})^{n} \cong j_{2\#}C(T^{2}_{\theta})^{n}$ with matrices in $GL_{n}(C(T^{2}_{\theta}))$ in the manner discussed in Section 2. Thus it is useful to have data concerning the group $\pi_{0}(GL_{n}(C(T^{2}_{\theta})))$ of path-components of $GL_{n}(C(T^{2}_{\theta}))$.

Following standard practice, we take the image in $GL_{n}(A)$ of an invertible element $x \in A$ to be the matrix $\begin{pmatrix} x & 0 \\ 0& 1_{n-1} \end{pmatrix}$.

\begin{prop} 
\label{UVgenerators}
The group $\pi_{0}(GL_{n}(C(T^{2}_{\theta})))\cong \mathbb{Z \times Z}$ is generated by the path-classes of the images of $U_{1}$ and $U_{2}$ in $GL_{n}(C(T_{\theta}^2))$. The natural map $\pi_{0}(GL_{n}(C(T^{2}_{\theta}))) \rightarrow K_{1}(C(T^{2}_{\theta}))$ is an isomorphism.
\end{prop}
\begin{proof} 
Since $V(C(T^{2}_{\theta}))$ and $V(C(T) \otimes C(T^{2}_{\theta}))$ are both cancellative (by Theorem 7.1 of \cite{projmodules}, if $\theta$ is irrational, and by the fact that $C(T^{2}_{\theta})$ is strongly Morita equivalent to $C(T^{2})$ if $\theta$ is rational), the natural map $\pi_{0}(GL_{n}(C(T^{2}_{\theta})))\rightarrow K_{1}(C(T^{2}_{\theta}))$ is injective by Theorem 8.4 of \cite{projmodules}, and thus an isomorphism, by Corollary 2.5 of \cite{PV}.
\end{proof}

This result generalizes to a noncommutative $m$-torus in the case that at least one entry of $\theta$ irrational:

\begin{theorem}
\label{rieffelgenerators}
(Theorem 3.3 of \cite{homotopy}). Let  $C(T^{m}_{\theta})$ be a noncommutative $m$-torus with not all entries of $\theta$ rational. Then
 \[\pi_{k}(GL_{n}(C(T^{m}_{\theta}))) \cong \left\{  \begin{array}{ll}
K_{1}(C(T^{m}_{\theta})) & \mbox{for k even}\\
K_{0}(C(T^{m}_{\theta})) & \mbox{for k odd}\end{array}  \right.
\cong \mathbb{Z}^{2^{m-1}}\] for all $k\geq 0$, $n \geq 1$.
\end{theorem}
The isomorphism is given by composing the natural map \[\pi_{k}(GL_{n}(C(T^{m}_{\theta})))\rightarrow  \pi_{k}(GL_{\infty}(C(T^{m}_{\theta})))\] induced by the usual embedding of $GL_{n}(C(T^{m}_{\theta}))$ into $GL_{\infty}(C(T^{m}_{\theta}))$ with Bott periodicity.

Notice that Theorem \ref{rieffelgenerators} is false in the commutative case $\theta = 0$ for $m \geq 3$, since $\pi_{0}(GL_{1}(C(T^{m}_{0}))) \cong [T^{m}, U(1)] \cong H^{1}(T^{m} ; \mathbb{Z}) \cong \mathbb{Z}^{m}$. In another direction, one can show using Postnikov approximation \cite{botttu} that $\pi_{0}(GL_{2}(C(T^{4}_{0}))) \cong [T^{4}, U(2)] \cong \mathbb{Z}^{8} \oplus \mathbb{Z}_{2}$.

\begin{theorem} 
\label{Amodstrivial}
All right finitely-generated projective $C(S^{3}_{\theta})$-modules are free. The free $C(S^{3}_{\theta})$-modules $C(S^{3}_{\theta})^{n}$ are mutually non-isomorphic.
\end{theorem}
\begin{proof}
As indicated in Proposition \ref{3spherepullback}, the C*-algebra $C(S^{3}_{\theta}) \cong \mathcal{A}$ is given by the pullback \[ \xymatrix{ \mathcal{A}\ar[d]_{i_2} \ar[r]^{i_1} &  \mathcal{A}_{1}  \ar[d]^{j_1} \\  \ \mathcal{A}_{2} \ar[r]_{j_2} & \mathcal{ B},}\] 
where the maps $i_{k}$ and $j_{k}$ are given by restriction.

All finitely-generated projective modules over $\mathcal{A}_{1}$ or $\mathcal{A}_{2}$ are free by Proposition \ref{filledmodstrivial}. Thus, by Theorem \ref{miln2}, and Propositions \ref{invtech1} and \ref{invdimension}, each finitely-generated projective $\mathcal{A}$-module is isomorphic to some $M(\mathcal{A}_{1}^{n}, \mathcal{A}_{2}^{n}, h)$, where $hj_{1\#}\mathcal{A}_{1}^{n}=j_{2\#}\mathcal{A}_{2}^{n}$.

As explained at the end of Section 2, we may use the standard ordered basis of $\mathcal{A}^{n}_{k}$ to identify automorphisms of $\mathcal{A}^{n}_{k}$ with matrices in $GL_{n}(\mathcal{A}_{k})$, and to identify isomorphisms  $j_{1\#}\mathcal{A}_{1}^{n} \cong j_{2\#}\mathcal{A}_{2}^{n}$ with matrices in $GL_{n}(\mathcal{B})$. As in Section 2, we write $id$ for the isomorphism $e_{1,i} \otimes_{\mathcal{A}_{1}} 1 \mapsto e_{2, i} \otimes_ {\mathcal{A}_{2}} 1$.

By Proposition \ref{UVgenerators}, there must be integers $k$ and $l$ such that $h$ is path-connected in $Iso_{\mathcal{B}}(j_{1\#}\mathcal{A}_{1}^{n}, j_{2\#}\mathcal{A}_{2}^{n}) \cong GL_{n}(\mathcal{B})$ to the isomorphism given by the image of $U_{1}^{k}U_{2}^{l}\in C(T^{2}_{\theta})$ in $GL_{n}(\mathcal{B})$. But this last isomorphism is just $j_{2*}(u_{1}^{k})\circ id \circ  j_{1*}((u_{2}^{-l})^{-1})$, where $u_{1}$ is the automorphism of $\mathcal{A}^{n}_{2}$ corresponding to the image of the constant function $t \mapsto U_{1}$ in $GL_{n}(\mathcal{A}_{2})$, and $u_{2}$ is the automorphism of $\mathcal{A}^{n}_{1}$ corresponding to the image of the constant function $t \mapsto U_{2}$ in $GL_{n}(\mathcal{A}_{1})$. So by Proposition \ref{clutching} and Corollary \ref{milncor},
\begin{gather*} M(\mathcal{A}_{1}^{n}, \mathcal{A}_{2}^{n}, h) \cong M(\mathcal{A}_{1}^{n}, \mathcal{A}_{2}^{n},  U_{1}^{k}U_{2}^{l}) \\ \cong M(\mathcal{A}_{1}^{n}, \mathcal{A}_{2}^{n}, id) \cong \mathcal{A}^{n} \cong C(S^{3}_{\theta})^{n}.\end{gather*} By Proposition \ref{invdimension}, the free $C(S^{3}_{\theta})$-modules $C(S^{3}_{\theta})^{n}$ are mutually non-isomorphic.

\end{proof}
\section{Finitely-Generated Projective Modules over $C(S^{4}_{\theta})$}
Throughout this section $\theta$ is an arbitrary real number unless otherwise specified. We classify and construct all finitely-generated projective $C(S^{4}_{\theta})$-modules up to isomorphism. Since $C(S^{4}_{\theta})$ is a pullback over $C(S^{3}_{\theta})$ (Proposition \ref{S4thetapullback}), we are able to use homotopy-theoretic data concerning $Iso_{C(S^{3}_{\theta})}(j_{1\#}C(D^{4}_{\theta})^{n}, j_{2\#}C(D^{4}_{\theta})^{n})\cong GL_{n}(C(S^{3}_{\theta}))$ to construct modules over $C(S^{4}_{\theta})$ by clutching. 

\begin{prop}
\label{D4trivial}
All finitely-generated projective $C(D^{2m}_{\theta})$-modules are free.
\end{prop}
\begin{proof}
The C*-algebra $C(D^{2m}_{\theta})$ deformation retracts onto its subalgebra of scalar multiples of the identity element. Explicitly, a retraction is given by the map \begin{gather*} F_{t}(w_{k})=(1-t)w_{k} \\F_{t}(y)=\sqrt{1-(1-t)^{2}(w_{1}w_{1}^{*}+...+ w_{m}w_{m}^{*})}
\end{gather*}
where the $w_{k}$ and $y$ are the generators of $C(D^{2m}_{\theta})$ given in definition \ref{evenball}.
\end{proof}

We need homotopy-theoretic data concerning the groups $GL_{n}(C(S^{3}_{\theta}))$. We first remind the reader of some standard constructions that can be found, for example, in \cite{blackadar}:

Let $A$ be a (not necessarily unital) C*-algebra. The group $GL_{n}(A)$ is defined as $\{ x \in GL_{n}(A^{\sim}) : x \equiv 1_{n}\mod M_{n}(A) \}$, where $A^{\sim}$ is the unitisation of $A$. This definition agrees with the usual notion of $GL_{n}(A)$ when $A$ is unital, but still makes sense when $A$ is nonunital. This allows nonunital C*-algebras to be incorporated into the $K$-theory.

Recall also that the suspension of $A$ is the nonunital C*-algebra $SA:=C_{0}((0,1), A)$. Bott periodicity gives us that
\[ K_{1}(A) \cong K_{0}(SA), \quad K_{0}(A) \cong K_{1}(SA). \]
The isomorphism $K_{0}(A) \cong K_{1}(SA)$ is given by the Bott map \[[e]-[1_{n}] \mapsto [(\exp(2 \pi i t)e+ 1_{n}-e)(\exp(-2 \pi i t)1_{n})],\] where $e$ is an idempotent in $M_{n}(A)$.

Thus, If $\theta$ is rational, then by strong Morita equivalence we have that \[K_{0}(SC(T^{2}_{\theta})) \cong K_{1}(C(T^{2}_{\theta})) \cong K^{1}(T^{2}) \cong \mathbb{Z \oplus Z}\] and \[K_{1}(SC(T^{2}_{\theta})) \cong K_{0}(C(T^{2}_{\theta})) \cong K^{0}(T^{2}) \cong \mathbb{Z \oplus Z}.\] 

If $\theta$ is irrational, then by corollary 2.5 of \cite{PV}, we again have that $K_{0}(SC(T^{2}_{\theta})) \cong \mathbb{Z \oplus Z}$ and, moreover, $K_{1}(SC(T^{2}_{\theta})) \cong \mathbb{Z \oplus Z}$ is generated by \[W:=\exp(2 \pi it)1 \in GL_{1}(SC(T^{2}_{\theta}))\] and \[X:=\exp(2 \pi it)p+1-p \in GL_{1}(SC(T^{2}_{\theta})),\] where $p \in C(T^{2}_{\theta})$ is a Rieffel projection of trace $\theta \mod 1$.

We abuse notation by using the symbols $W$ and $X$ to also refer to the images of $W$ and $X$ in $GL_{n}(SC(T^{2}_{\theta}))$ and $K_{1}(SC(T^{2}_{\theta}))$.

\begin{prop}
\label{suspensiongen}
Let $\theta$ be irrational. Then, for any $n \geq 1$, the group $\pi_{0}(GL_{n}(SC(T^{2}_{\theta}))) \cong \mathbb{Z \times Z}$ is generated by $W$ and $X$.
\end{prop}
\begin{proof}
Let $A$ be a unital C*-algebra. Regarding an element of $GL_{n}(SA)$ as a base-pointed loop (with base-point the identity of A) through $GL_{n}(A)$, we have an isomorphism $\pi_{0}(GL_{n}(SA)) \cong \pi_{1}(GL_{n}(A))$.

Thus we have a commutative diagram \[ \xymatrix{\pi_{0}(GL_{n}(SC(T^{2}_{\theta})))
 \ar[d] \ar[r]^{\cong} & \pi_{1}(GL_{n}(C(T^{2}_{\theta}))) \ar[d] \\ K_{1}(SC(T^{2}_{\theta})) \ar[r]_{\cong} & K_{0}(C(T^{2}_{\theta})).}\]
The map $\pi_{1}(GL_{n}(C(T^{2}_{\theta}))) \rightarrow K_{0}(C(T^{2}_{\theta}))$ is an isomorphism by Theorem \ref{rieffelgenerators}. Thus the result follows from the remarks directly preceeding the proposition.
\end{proof}

Now returning to the analysis of the groups $\pi_{0}(GL_{n}(C(S^{3}_{\theta})))$, we first observe that it follows immediately from Proposition \ref{3spherepullback} that there is a short exact sequence \[ 0 \rightarrow SC(T_{\theta}^{2}) \xrightarrow{i} C(S^{3}_{\theta}) \xrightarrow{q} C^{*}(U_{2}) \oplus C^{*}(U_{1}) \rightarrow 0. \]

\begin{lemma}
\label{surjective}
The induced map $ i_{*}: K_{1}(SC(T^{2}_{\theta})) \rightarrow K_{1}(C(S^{3}_{\theta}))$ is surjective.
\end{lemma}
\begin{proof} 
The the standard six-term exact sequence of Banach algebra $K$-theory (theorem 9.3.1 of \cite{blackadar}) gives us \[ \xymatrix{ K_{0}(SC(T^{2}_{\theta}))\ar[r]^-{i_{*}}& K_{0}(C(S^{3}_{\theta})) \ar[r]^-{q_{*}}  & K_{0}(C^{*}(U_{2}))\oplus K_{0}(C^{*}(U_{1})) \ar[d]^-{\partial} \\
K_{1}(C^{*}(U_{2}))\oplus K_{1}(C^{*}(U_{1})) \ar[u]^-{\partial} &  K_{1}(C(S^{3}_{\theta})) \ar[l]_-{q_{*}} & K_{1}(SC(T^{2}_{\theta})) \ar[l]_-{i_{*}}} \]

or,  \[ \xymatrix{\mathbb{Z \oplus Z} \ar[r]^-{i_{*}} & \mathbb{Z} \ar[r]^-{q_{*}}  & \mathbb{Z \oplus Z} \ar[d]^-{\partial} \\
\mathbb{Z \oplus Z} \ar[u]^-{\partial} &  \mathbb{Z} \ar[l]_-{q_{*}} & \mathbb{Z \oplus Z}. \ar[l]_-{i_{*}}}\]
\newline

The above values of the $K$-groups, the condition of exactness, and the fact that the map $q_{*}: K_{0}(C(S^{3}_{\theta})) \rightarrow K_{0}(C^{*}(U_{2}))\oplus K_{0}(C^{*}(U_{1}))$ is injective (established below), jointly suffice to completely determine the images and kernels (at least up to isomorphism as abstract discrete groups) of all of the other maps in the six-term sequence.

Indeed, from the Grothendeick construction, the map \[q_{*}:K_{0}(C(S^{3}_{\theta})) \rightarrow K_{0}(C^{*}(U_{2}))\oplus K_{0}(C^{*}(U_{1}))\] is induced by the semigroup  homomorphism \[ q_{*}:V(C(S^{3}_{\theta})) \rightarrow V(C^{*}(U_{2}))\oplus V(C^{*}(U_{1})).\] 

The algebras $C(S^{3}_{\theta})$, $C^{*}(U_{1})$, and $C^{*}(U_{2})$ all satisfy the invariance of dimension property, and every finitely-generated projective module over one of these algebras is free. Thus, identifying the rank of a free module over one of these algebras with the isomorphism class of that free module, it is clear both that the map $q_{*}:V(C(S^{3}_{\theta})) \rightarrow V(C^{*}(U_{2}))\oplus V(C^{*}(U_{1}))$ is simply $[n] \mapsto ([n], [n])$, for $n$ any non-negative integer, and, moreover, it's extension $q_{*}:K_{0}(C(S^{3}_{\theta})) \rightarrow K_{0}(C^{*}(U_{2}))\oplus K_{0}(C^{*}(U_{1}))$ is $[k] \mapsto ([k],[k])$, for $k$ any integer. Therefore, by exactness, the map $i_{*}:K_{0}(SC(T^{2}_{\theta})) \rightarrow K_{0}(C(S^{3}_{\theta}))$ is the zero map, and the index map $\partial : K_{1}(C^{*}(U_{2}))\oplus K_{1}(C^{*}(U_{1})) \rightarrow K_{0}(SC(T^{2}_{\theta}))$ is surjective. But given exactness and the values (as abstract discrete groups) of the $K$-groups in the six-term sequence, the index map cannot be surjective unless the map $q_{*}: K_{1}(C(S^{3}_{\theta})) \rightarrow   K_{1}(C^{*}(U_{2}))\oplus K_{1}(C^{*}(U_{1}))$ is the zero map. So $ i_{*}: K_{1}(SC(T^{2}_{\theta})) \rightarrow K_{1}(C(S^{3}_{\theta}))$ is surjective.
\end{proof}

We continue to abuse notation by also referring to the images of $W$ and $X$ in $GL_{n}(C(S^{3}_{\theta}))$ and $K_{1}(C(S^{3}_{\theta}))$ by $W$ and $X$.

\begin{lemma} 
\label{pathtoidentity}
The invertible $W$  is path-connected through $GL_{1}(C(S^{3}_{\theta}))$ to the identity element.
\end{lemma}
\begin{proof}
The map $t \mapsto 2\pi it$ is in $\mathcal{A} \cong C(S^{3}_{\theta})$. So $W$ is path-connected in $GL_{1}(\mathcal{A})$ to the constant function $1 \in \mathcal{A}$, since $W$ is the exponential of a function in $\mathcal{A}$.
\end{proof}

\begin{lemma} 
\label{Agenerator}
If $\theta$ is irrational, then $X$ generates $K_{1}(C(S^{3}_{\theta})) \cong \mathbb{Z}$.
\end{lemma}
\begin{proof}
 Immediate by Lemmas \ref{surjective} and \ref{pathtoidentity}.
\end{proof}

\begin{theorem}
\label{naturalmap}
If $\theta$ is irrational, then $X$ generates $\pi_{0}(GL_{1}(C(S^{3}_{\theta})))$. The natural map $\pi_{0}(GL_{n}(C(S^{3}_{\theta}))) \rightarrow K_{1}(C(S^{3}_{\theta}))$ is an isomorphism for all $n \geq 1$.
\end{theorem}
\begin{proof}
For each $n \geq 1$, we have a commutative diagram with exact rows:
 \[ \xymatrix{\pi_{0}(GL_{n}(SC(T^{2}_{\theta})))\ar[r]^-{i_{*}} \ar[d]^{\cong} & \pi_{0}(GL_{n}(C(S^{3}_{\theta}))) \ar[r]^-{q_{*}} \ar[d] & \pi_{0}(GL_{n}(C^{*}(U_{2})))\oplus \pi_{0}(GL_{n}(C^{*}(U_{1}))) \ar[d]^{\cong} \\
\ K_{1}(SC(T^{2}_{\theta}))\ar[r]^-{i_{*}}& K_{1}(C(S^{3}_{\theta})) \ar[r]^-{q_{*}}  & K_{1}(C^{*}(U_{2}))\oplus K_{1}(C^{*}(U_{1})). }\]

The natural map \[ \pi_{0}(GL_{n}(C^{*}(U_{2})))\oplus \pi_{0}(GL_{n}(C^{*}(U_{1}))) \rightarrow  K_{1}(C^{*}(U_{2}))\oplus K_{1}(C^{*}(U_{1}))\] is an isomorphism, since $C^{*}(U_{1})$ and $C^{*}(U_{2})$ are isomorphic to the circle algebra $C(S^{1})$.

The natural map $\pi_{0}(GL_{n}(SC(T^{2}_{\theta}))) \rightarrow K_{1}(SC(T^{2}_{\theta}))$ is an isomorphism by Proposition \ref{suspensiongen}.

Suppose now that there is an element $Y \in GL_{n}(C(S^{3}_{\theta}))$ that is not homotopic to the image of $X^{k}=(\exp(2 \pi it)p+1-p)^{k}$ in $GL_{n}(C(S^{3}_{\theta}))$ for some $k$. Then, since $X$ and $W$ generate $\pi_{0}(GL_{n}(SC(T^{2}_{\theta})))$, we conclude by Lemma \ref{pathtoidentity} that $Y$ is not in the image of $\pi_{0}(GL_{n}(SC(T^{2}_{\theta})))$ in $ GL_{n}(C(S^{3}_{\theta}))$, and so $q_{*}(Y)$ is not trivial in $\pi_{0}(GL_{n}(C^{*}(U_{2})))\oplus \pi_{0}(GL_{n}(C^{*}(U_{1}))).$ But $q_{*}(Y)$ is trivial in $K_{1}(C^{*}(U_{2}))\oplus K_{1}(C^{*}(U_{1}))$, by Lemma \ref{surjective}, and so we have a contradiction by the commutativity of the diagram. 
\end{proof}

We are now in position to classify and construct all finitely-generated projective $C(S^{4}_{\theta})$-modules up to isomorphism. Let $\theta$ be irrational. Following the convention of Section 2, we write $M(C(D^{4}_{\theta})^{n}, C(D^{4}_{\theta})^{n}, X^{s})$ as short-hand for $M(C(D^{4}_{\theta})^{n},C(D^{4}_{\theta})^{n}, h)$, where $h$ is the isomorphism $j_{1\#}C(D^{4}_{\theta})^{n} \cong j_{2\#}C(D^{4}_{\theta})^{n}$ that corresponds to the image of  $X^{s}=(\exp(2\pi it)p+1-p)^{s}$ in $ GL_{n}(C(S^{3}_{\theta}))$ under our identification \[Iso_{C(S^{3}_{\theta})}(j_{1\#}C(D^{4}_{\theta})^{n}, j_{2\#}C(D^{4}_{\theta})^{n}) \cong GL_{n}(C(S^{3}_{\theta})),\] and where the $j_{k}$'s are as in Proposition \ref{S4thetapullback}.

\begin{definition}
\label{s4modules}
Let $\theta$ be irrational. We define $N(n,s)$ to be the finitely-generated projective $C(S^{4}_{\theta})$-module $M(C(D^{4}_{\theta})^{n}, C(D^{4}_{\theta})^{n}, X^{s})$.
\end{definition}

\begin{theorem}
\label{maintheorem}
Let $\theta$ be irrational. Then the map $[M(n,s)] \rightarrow (n,s)$ gives a semigroup isomorphism \[ V(C(S^{4}_{\theta}))\cong \{0\} \cup (\mathbb{N}\times K_{1}(C(S^{3}_{\theta}))) \cong \{0\} \cup (\mathbb{N} \times \mathbb{Z}).\]
\end{theorem}
\begin{proof}
We denote the two copies of $C(D^{4}_{\theta})$ by $A_{1}$ and $A_{2}$, as in Proposition \ref{S4thetapullback}. We have observed in Proposition \ref{D4trivial} that all finitely-generated projective $A_{k}$-modules are free. Also the $A_{k}$ satisfy the invariance of dimension property, since they admit tracial states. So by Theorem \ref{miln2} and Proposition \ref{invtech1}, every finitely-generated projective $C(S^{4}_{\theta})$-module is isomorphic to one of the form $M(A_{1}^{n}, A_{2}^{n}, h)$. But since the free $C(S^{3}_{\theta})$-modules $j_{1\#}A_{1}^{n}$ and $j_{2\#}A_{2}^{n}$ are isomorphic, we may follow Section 2 and identify the $j_{k\#}A_{k}^{n}$'s with $C(S^{3}_{\theta})^{n}$, and then identify the isomorphism $h$ with an element of $GL_{n}(C(S^{3}_{\theta}))$. But, by Theorems \ref{Agenerator} and \ref{naturalmap}, each $h$ is path-connected in $ Iso_{C(S^{3}_{\theta})}(j_{1\#}A_{1}^{n}, j_{2\#}A_{2}^{n})\cong GL_{n}(C(S^{3}_{\theta}))$ to the isomorphism corresponding to $\begin{pmatrix} X^{s} & 0 \\ 0 & 1_{n-1}\end{pmatrix}$, for some unique integer $s$. Thus each finitely-generated projective $C(S^{4}_{\theta})$-module is isomorphic to some module $N(n, s)$ by Proposition \ref{clutching}.

Suppose now that $N(n_{1},s_{1})$ is isomorphic to $N(n_{2}, s_{2})$. Then $n_{1} = n_{2}$, by Proposition \ref{invtech2}.

Now suppose that $N(n,s_{1})$ is isomorphic to $N(n,s_{2})$. Then, applying Corollary \ref{milncor}, and writing matrices for the corresponding elements of $Iso_{C(S^{3}_{\theta})}(j_{1\#}A_{1}^{n}, j_{2\#}A_{2}^{n})$, we have that \[ \begin{pmatrix} X^{s_{1}} & 0 \\ 0 & 1_{n-1}\end{pmatrix}=j_{2*}(g_{2})\begin{pmatrix} X^{s_{2}} & 0 \\ 0 & 1_{n-1}\end{pmatrix}j_{1*}(g_{1}^{-1}),\] for some $g_{k} \in Aut_{A_{k}}(A_{k}^{n}) \cong GL_{n}(A_{k})$. But $GL_{n}(A_{k})$ is path-connected, so $j_{2*}(g_{2})$ is path-connected in $Aut_{C(S^{3}_{\theta})}(j_{2\#}A_{2}^{n})$ to the identity automorphism of $A_{2}^{n}$, and $j_{1*}(g_{1}^{-1})$ is path-connected in $Aut_{C(S^{3}_{\theta})}(j_{1\#}A_{1}^{n})$  to the identity automorphism of $A_{1}^{n}$. But this means that \[\begin{pmatrix} X^{s_{1}} & 0 \\ 0 & 1_{n-1}\end{pmatrix}=
j_{2*}(g_{2})\begin{pmatrix} X^{s_{2}} & 0 \\ 0 & 1_{n-1}\end{pmatrix}j_{1*}(g_{1}^{-1})\sim_{h} \begin{pmatrix} X^{s_{2}} & 0 \\ 0 & 1_{n-1}\end{pmatrix},\]which is possible only if $s_{1}=s_{2}$, by Lemma \ref{Agenerator} and Theorem \ref{naturalmap}.

Finally, since $\begin {pmatrix} X^{s_{1}} & 0\\ 0 & X^{s_{2}} \end{pmatrix}$ is homotopic through $GL_{n_{1}+n_{2}}(C(S^{3}_{\theta}))$ to $\begin{pmatrix} X^{s_{1}}X^{s_{2}} & 0 \\ 0 & 1_{n_{1}+n_{2}-1}\end{pmatrix}$ , we see that
\begin{equation*}
\begin{split}
N(n_{1}, s_{1}) \oplus N(n_{2}, s_{2}) & \cong N(n_{1}+n_{2}, s_{1}+s_{2})\\
& \cong N(1,s_{1}+s_{2})\oplus C(S^{4}_{\theta})^{n_{1}+n_{2}-1}.
\end{split}
\end{equation*}
\end{proof}

We note that since $V(C(S^{2}))$ is also $\{0\} \cup (\mathbb{N} \times \mathbb{Z})$, we have that $V(C(S^{4}_{\theta}))$ and $V(C(S^{2}))$ are isomorphic. Indeed, the complex vector bundles over $S^{2}$ are indexed (up to isomorphism) by their ranks and Chern numbers (integral of the first Chern class), so that every complex vector bundle over $S^{2}$ is isomorphic to a sum of a line bundle and a trivial bundle. The structural parallelism between $V(C(S^{4}_{\theta}))$ and $V(C(S^{2}))$ is not a mere formal coincidence, but rather is of a geometric nature as we are now in a position to show.

First recall that $C(D^{4}_{\theta}) \cong (C_{0}([0,1))\otimes C(S^{3}_{\theta}))^{\sim}$ by Corollary \ref{cone}.
Now consider the element $X=\exp(2\pi it)p+1-p$ in $C(\{0\}) \otimes C(S^{3}_{\theta})) \cong C(S^{3}_{\theta})$. Lift $p$ to the self-adjoint element $q=f\otimes p \in C_{0}([0,1))\otimes C(S^{3}_{\theta})\subset A_{1}= C(D^{4}_{\theta})$, where $f$ is the function $u \mapsto 1-u$ on $[0,1)$. The normal element $c=\exp(2\pi it)q+1-q$ is then a lift of $X$ to $A_{1}$. Since $p$ is a nontrivial projection, the spectrum $\sigma (X)$ is the full circle $S^{1}$ and  $\sigma (q)= [0,1]$. For each fixed $t$, the spectrum of $\exp(2 \pi it)q+1-q$ is the chord that connects the point $1$ to the point $\exp (2 \pi it)$, so $\sigma(c)=D^{2}$. 

We can similarly lift $p$ to a self-adjoint element $q^{\prime}$ in $A_{2}=C(D^{4}_{\theta})$ such that $\sigma (q^{\prime})=[0, 1]$, and define the normal element $c^{\prime}=\exp(2 \pi it)q^{\prime}+1-q^{\prime}$. As before,  the spectrum of $c^{\prime}$ is $D^{2}$.

From the above considerations of spectrum, the pullback $\mathcal P$ of $C^{*}(c,1)$ and $C^{*}(c^{\prime}, 1)$ over $C^{*}(X,1)$ that comes from restriction of the pullback diagram of Proposition \ref{S4thetapullback}, is a C*-subalgebra of $C(S^{4}_{\theta})$ that is  isomorphic to $C(S^{2})$.

\begin{cor}
\label{tensor}
Suppose that $\theta$ is irrational. Then \[N(n,s) \cong M(n,s)\otimes_{C(S^{2})} C(S^{4}_{\theta}),\] where $M(n,s)$ is the module of continuous sections of a rank-$n$ complex vector bundle over $S^{2}$ with Chern number $-s$, and the inclusion $C(S^{2}) \cong \mathcal{P} \hookrightarrow C(S^{4}_{\theta})$ is as described above.
\end{cor}
\begin{proof}

We form the $\mathcal P$-module $M(C^{*}(c,1)^{n}, C^{*}(c^{\prime}, 1)^{n}, X^{s})$, where we view $X^{s}\in GL_{n}(C^{*}(X,1))$ as in $Iso_{C^{*}(X, 1)}(j_{1\#}C^{*}(c,1)^{n}, j_{2\#}C^{*}(c^{\prime}, 1)^{n})$. Let $i:\mathcal P \hookrightarrow C(S^{4}_{\theta})$ be the inclusion map coming from the isomorphism of  Proposition \ref{S4thetapullback}. We first need to see that the $C(S^{4}_{\theta})$-module $i_{\#}M(C^{*}(c,1)^{n}, C^{*}(c^{\prime}, 1)^{n}, X^{s})$ is isomorphic to $N(n,s)$. 

Using $i$ now to also refer to the inclusions \[C^{*}(c,1) \hookrightarrow A_{1}, \quad C^{*}(c^{\prime}, 1) \hookrightarrow A_{2}, \quad C^{*}(X,1) \hookrightarrow C(S^{3}_{\theta}),\] it is clear that $i_{\#}C^{*}(c, 1)^{n} \cong A_{1}^{n}$ and $i_{\#}C^{*}(c^{\prime}, 1)^{n} \cong A_{2}^{n}$. Viewing $X^{s}$ as in $Iso_{C^{*}(X, 1)}(j_{1\#}C^{*}(c,1)^{n}, j_{2\#}C^{*}(c^{\prime}, 1)^{n})$, it is clear that $i_{*}(X^{s}) = X^{s} \otimes id$ identifies with $X^{s}$ viewed as in $Iso_{C(S^{3}_{\theta})}(j_{1\#}A_{1}^{n}, j_{2\#}A_{2}^{n})$, because $i_{\#}j_{1\#}C^{*}(c,1)^{n} \cong i_{\#}C^{*}(X,1)^{n}= C^{*}(X, 1)^{n} \otimes_{C^{*}(X,1)} C(S^{3}_{\theta}) \cong C(S^{3}_{\theta})^{n}$, and, similarly,  $i_{\#}j_{2\#}C^{*}(c^{\prime},1)^{n} \cong C(S^{3}_{\theta})^{n}$. 
Thus the module $i_{\#}M(C^{*}(c,1)^{n}, C^{*}(c^{\prime}, 1)^{n}, X^{s})$ is isomorphic to $N(n,s)$.

 But under the isomorphism $C^{*}(X,1) \cong C(S^{1})$, multiplication by $X^{s}$ corresponds to multiplication of functions in $C(S^{1})$ by the function $z=\exp(2 \pi it) \mapsto z^{s}=\exp(2 \pi ist)$. So $M(C^{*}(c,1)^{n}, C^{*}(c^{\prime}, 1)^{n}, X^{s})$ is isomorphic to the space of continuous sections of the  complex rank-$n$ vector bundle over $S^{2}$ that is formed by  by using the clutching function $z \mapsto \begin{pmatrix} z^{s} &  0 \\ 0 & 1_{n-1} \end{pmatrix}$ to glue a rank-$n$ trivial bundle over the northern hemisphere to a rank-$n$ trivial bundle over the southern hemisphere. (This vector bundle has Chern number equal to $-s$, if the equator $S^{1}$ inherits its orientation from the standard orientation of $S^{2}$ in the usual way.)

Thus, since each finitely-generated projective $C(S^{4}_{\theta})$-module is isomorphic to one of the $N(n,s)$, the inclusion $i:C(S^{2}) \cong \mathcal P \hookrightarrow C(S^{4})$ induces a semigroup isomorphism $i_{*}: V(C(S^{2})) \cong V(C(S^{4}_{\theta}))$ given by $M(C(D^{2})^{n}, C(D^{2})^{n},  z^{s}) \mapsto N(n,s)$.
\end{proof}

Heuristically speaking, Corollary \ref{tensor} states that the isomorphism classes of ``complex vector bundles" over the ``noncommutative space" $S^{4}_{\theta}$ correspond to the isomorphism classes of complex vector bundles over $S^{2}$ via pullback of bundles by a certain fixed quotient map $S^{4}_{\theta} \rightarrow S^{2}$.

Of course Corollary \ref{tensor} is false in the classical case $\theta=0$ (and also if $\theta$ is any rational number, as follows easily from Theorem \ref{rationalcase} below). This is because the pullback of any complex line bundle over $S^{2}$ by any quotient map $S^{4} \rightarrow S^{2}$ must be a trivial bundle, since there are no nontrivial complex line bundles over $S^{4}$. But then the pullback of any complex vector bundle over $S^{2}$ by such a map is a trivial vector bundle, since every complex vector bundle over $S^{2}$ is isomorphic to the direct sum of a line bundle and a trivial bundle.
\newline

We now consider the rational case $\theta=p/q$. We first observe that $C(S^{4}_{p/q})$ and $C(S^{4})$ are strongly Morita equivalent.

\begin{prop}
\label{moritaequiv}
Suppose that all entries of $\theta$ are rational. Then $C(S^{2m-1}_{\theta})$ is strongly Morita equivalent to $C(S^{2m-1})$ and $C(S^{2m}_{\theta})$ is strongly Morita equivalent to $C(S^{2m}).$
\end{prop}
\begin{proof}

Since all entires of $\theta$ are rational, we have that $C(T^{m}_{\theta})$ and $C(T^{m})$ are strongly Morita equivalent, by corollary 4.2 of \cite{elliottli}. Equivalently, we have that $C(T^{m}_{\theta})\otimes \mathbb{K}$ is isomorphic to $C(T^{m}) \otimes \mathbb{K}$, where $\mathbb{K}$ is the algebra of compact operators on a separable infinite-dimensional Hilbert space. Thus 
\begin{align*} 
 C(S^{2m-1}_{\theta}) \otimes \mathbb{K} & \cong \\
\{f \in C(S^{m-1}_{+}, C(T^{m}_{\theta})) : f(S(i)) \subset C(T^{m}_{\theta})(i) \text{ all } i \leq m\} \otimes \mathbb{K} & \cong \\
\{f \in C(S^{m-1}_{+}, C(T^{m}_{\theta}) \otimes \mathbb{K}) : f(S(i)) \subset C(T^{m}_{\theta})(i) \otimes \mathbb{K}  \text{ all } i \leq m\} & \cong \\ \{f \in C(S^{m-1}_{+}, C(T^{m}) \otimes \mathbb{K}) : f(S(i)) \subset C(T^{m})(i) \otimes \mathbb{K}  \text{ all } i \leq m\} & \cong \\ \{f \in C(S^{m-1}_{+}, C(T^{m})) : f(S(i)) \subset C(T^{m})(i) \text{ all } i \leq m\} \otimes \mathbb{K} & \cong \\ C(S^{2m-1}) \otimes \mathbb{K},
\end{align*}
The first and last isomorphisms here are by Theorem \ref{oddfield}.  The third isomorphism follows from the fact that the isomorphism $C(T^{m}_{\theta}) \otimes \mathbb{K} \cong C(T^{n}) \otimes \mathbb{K}$ restricts to an isomorphism $C(T^{m}_{\theta})(i) \otimes \mathbb{K} \cong C(T^{m})(i)\otimes \mathbb{K}$.
Thus $C(S^{2m-1}_{\theta})$ and $C(S^{2m-1})$ are strongly Morita equivalent.

Similarly,
\begin{align*} 
 C(S^{2m}_{\theta}) \otimes \mathbb{K} & \cong \\
\{f \in C(S^{m}_{+, m+1}, C(T^{m}_{\theta})) : f(S(i)) \subset C(T^{m}_{\theta})(i) \text{ all } i \leq m\} \otimes \mathbb{K} & \cong \\
\{f \in C(S^{m}_{+, m+1}, C(T^{m}_{\theta}) \otimes \mathbb{K}) : f(S(i)) \subset C(T^{m}_{\theta})(i) \otimes \mathbb{K}  \text{ all } i \leq m\} & \cong \\ \{f \in C(S^{m}_{+, m+1}, C(T^{m}) \otimes \mathbb{K}) : f(S(i)) \subset C(T^{m})(i) \otimes \mathbb{K}  \text{ all } i \leq m\} & \cong \\ \{f \in C(S^{m}_{+, m+1}, C(T^{m})) : f(S(i)) \subset C(T^{m})(i) \text{ all } i \leq m\} \otimes \mathbb{K} & \cong \\ C(S^{2m}) \otimes \mathbb{K},
\end{align*}

 so $C(S^{2m}_{\theta})$ and $C(S^{2m})$ are strongly Morita equivalent.

\end{proof}

\begin{theorem}
\label{rationalcase}
Let $\theta=p/q$ be rational. Then there is a semigroup isomorphism \[ V(C(S^{4}_{p/q}))\cong \{0, 1\}\cup (\mathbb{N}\times K_{1}(C(S^{3}_{\theta}))) \cong \{0, 1\} \cup ((\mathbb{N} \backslash \{1\}) \times \mathbb{Z})).\]
Every finitely-generated projective $C(S^{4}_{p/q})$-module is either a rank-1 free module or is isomorphic to a module of the form \[M(C(D^{4}_{p/q})^{n}, C(D^{4}_{p/q})^{n}, \begin{pmatrix} z_{2} & z_{1} \\ -\rho z_{1}^{*} & z_{2}^{*} \end{pmatrix}^{s}),\] for some $n\geq 2$ (here identifying  $\begin{pmatrix} z_{2} & z_{1} \\ -\rho z_{1}^{*} & z_{2}^{*} \end{pmatrix}^{s})$ with its image in $Iso_{C(S^{3}_{p/q})}(j_{1\#}C(D^{4}_{p/q})^{n}, j_{2\#}C(D^{4}_{p/q})^{n}) \cong GL_{n}(C(S^{3}_{p/q}))$). Thus the isomorphism $ V(C(S^{4}_{p/q}))\cong \{0, 1\} \cup ((\mathbb{N} \backslash \{1\}) \times \mathbb{Z}))$ is given by  \[ [M(C(D^{4}_{p/q})^{n}, C(D^{4}_{p/q})^{n}, \begin{pmatrix} z_{2} & z_{1} \\ -\rho z_{1}^{*} & z_{2}^{*} \end{pmatrix}^{s})] \mapsto (n,s).\] 
\end{theorem}
\begin{proof}
It follows from the proof of Theorem 5.5 of \cite{natsume} that $\begin{pmatrix} z_{2} & z_{1} \\ -\rho z_{1}^{*} & z_{2}^{*} \end{pmatrix}$ generates $K_{1}(C(S^{4}_{p/q}))$ (see the discussion on page 28). The semigroup $V(C(T)\otimes C(S^{3}_{p/q}))$ satisfies cancellation since $C(T)\otimes C(S^{3}_{p/q})$ is strongly Morita equivalent to $C(T\times S^{3})$. Thus the image of $\begin{pmatrix} z_{2} & z_{1} \\ -\rho z_{1}^{*} & z_{2}^{*} \end{pmatrix}$ in $GL_{n}(C(S^{3}_{p/q}))$ must generate $\pi_{0}(GL_{n}(C(S^{3}_{p/q})))$, for all $n\geq 2$, by Theorem 8.4 of \cite{projmodules}.

 Every finitely-generated projective $C(S^{4}_{p/q})$-module (up to isomorphism) comes from gluing together two rank-$n$ free $C(D^{4}_{p/q})$-modules by some element of  $GL_{n}(C(S^{3}_{p/q}))$ for some $n$. But  $V(C(S^{4}_{p/q})) \cong V(C(S^{4}))$ since $C(S^{4}_{p/q})$ is strongly Morita equivalent to $C(S^{4})$. From these last two observations, and the known structure of the semigroup $V(C(S^{4}))$,  we see that $\pi_{0}(GL_{1}(C(S^{3}_{p/q}))) \ncong 0$ is impossible. 

Thus every finitely-generated projective $C(S^{4}_{p/q})$-module is isomorphic to either a rank-1 free $C(S^{4}_{p/q})$-module or to a module of the form $M(C(D^{4}_{p/q})^{n}, C(D^{4}_{p/q})^{n}, \begin{pmatrix} z_{2} & z_{1} \\ -\rho z_{1}^{*} & z_{2}^{*} \end{pmatrix}^{s}),$ for some $n\geq 2$. 

\end{proof}

Thus, for both the rational and irrational cases, we have completely characterized the set of isomorphism classes of finitely-generated projective $C(S^{4}_{\theta})$-modules as a semigroup. We immediately obtain:

\begin{cor}
\label{cancellationthm}
Let $\theta$ be any real number. The algebra $C(S^{4}_{\theta})$ is cancellative.
\end{cor}

We relate the $N(n,s)$ to the $C(S^{4}_{\theta})$-modules already appearing in the literature.

First suppose that in the context of Section 2, we are given an $A$-module $PA^{N}$, where both $A_{k}$ satisfy the invariance of dimension property, and $P$ is an idempotent in $M_{n}(A)$. Suppose that $i_{k\#}PA^{N} \cong A_{k}^{n}$, and $j_{k\#}i_{k\#}PA^{N} \cong B^{n}$. We may then construct an $h \in GL_{n}(B) \cong Aut_{B}(B^{n}) \cong Iso_{B}(j_{1\#}i_{1\#}PA^{N}, j_{2\#}i_{2\#}PA^{N})$  so that $PA^{N}$ is isomorphic to $M(A_{1}^{n}, A_{2}^{n}, h)$, as follows:
We have a commutative diagram  \[ \xymatrix{PA^{N} \ar[d] \ar[r] & i_{1}(PA^{N})  \ar[d] \ar[r]^{\psi_{1}}_{\cong} &  A_{1}^{n} \ar[d] \\ i_{2}(PA^{N}) \ar[r] \ar[d]_{\psi_{2}}^{\cong} & j_{2}i_{2}(PA^{N})= j_{1}i_{1}(PA^{N}) \ar[d]_{j_{2}(\psi_{2})}^{\cong} \ar[r]^-{j_{1}(\psi_{1})}_-{\cong} & B^{n}
\\ A_{2}^{n} \ar[r] & B^{n}
.}\]
Here we identify the module $i_{k\#}PA^{N}$  with $i_{k}(PA^{N})$, and also identify the module $j_{k\#}i_{k\#}PA^{N}$ with $j_{k}i_{k}(PA^{N})$, and the canonical map $f_{PA^{N}}$ with the identity map $j_{1}i_{1}(PA^{N}) \rightarrow j_{2}i_{2}(PA^{N})=j_{1}i_{1}(PA^{N})$.

We also identify $A^{n}_{k} \otimes_{A_{k}} B$ with $B^{n}$. Commutativity of the diagram then forces the maps $A^{n}_{k}\rightarrow B^{n}$ to be the maps that take the $i$-th member of the standard ordered basis of $A_{k}^{n}$ to the $i$-th member of the standard ordered basis of $B^{n}$. By Corollary \ref{milncor}, if we take $h$ to be the automorphism $j_{2}(\psi_{2}) j_{1}(\psi_{1})^{-1}$, then $PA^{N} \cong M(A_{1}^{n}, A_{2}^{n}, h)$.

The charge 1 instanton bundle $eC(S^{4}_{\theta})^{4}$  of Connes and Landi is given by the projection
 \[ e:= \frac{1}{2} \begin{pmatrix} 1+x & 0 & z_{2} & z_{1} \\ 0 & 1+x & -\rho z_{1}^{*} &  z_{2}^{*} \\ z_{2}^{*} & -\bar{\rho} z_{1} & 1-x & 0\\ z_{1}^{*} & z_{2} &  0 & 1-x\end{pmatrix},\] where $\rho=\exp(2 \pi i \theta )$, and the $z_{i}$ and $x$ are the generators of $C(S^{4}_{\theta})$ given in Definition \ref{evenspheres}.

Let $A=C(S^{4}_{\theta})$ and $B=C(S^{3}_{\theta})$. Let the $A_{k}$ and the $i_{k}$'s and $j_{k}$'s be as in Proposition \ref{S4thetapullback}. Since the $i_{k}(eC(S^{4}_{\theta})^{4})$ are free $A_{k}$-modules, and since the projections $i_{k}(e)$ have trace equal to 2, for any trace on $M_{4}(A_{k})$ that comes from any normalized tracial state on $A_{k}$, we know \emph{a priori} that $i_{k}(eC(S^{4}_{\theta})^{4})$ is isomorphic to $A_{k}^{2}$. Thus to explicitly find an $h$ so that $eC(S^{4}_{\theta})^{4} \cong M(A_{1}^{2}, A_{2}^{2}, h)$, it suffices to explicitly find trivializations $\psi_{k}: i_{k}(eC(S^{4}_{\theta})^{4}) \cong A_{k}^{2}$.

View elements of $A_{k}^{2}$ as column vectors and let $\psi_{1}$ be multiplication by \[\begin{pmatrix} 1 & 0 & (1+y)^{-1}w_{2} & (1+y)^{-1}w_{1} \\ 0 & 1 & -\rho (1+y)^{-1}w_{1}^{*} &  (1+y)^{-1}w_{2}^{*}  \end{pmatrix}. \]

 It is easy to check that multiplication by \[ \frac{1}{2} \begin{pmatrix} 1+y & 0 \\ 0 & 1+y \\ w_{2}^{*} & - \bar{\rho}w_{1} \\ w_{1}^{*} &  w_{2} \end{pmatrix} \] gives the inverse to $\psi_{1}$.

We let $\psi_{2}$ be multiplication by \[\begin{pmatrix} (1+y)^{-1}w_{2}^{*}  & -\bar{\rho}(1+y)^{-1}w_{1}  & 1 & 0  \\ (1+y)^{-1}w_{1}^{*} & (1+y)^{-1}w_{2} & 0 &  1 \end{pmatrix}. \]

Multiplication by \[ \frac{1}{2} \begin{pmatrix} w_{2}& w_{1} \\ -\rho w_{1}^{*} & w_{2}^{*}\\ 1+y & 0 \\ 0 &  1+y \end{pmatrix} \] gives the inverse to $\psi_{2}$.

Direct calculation gives us \[h=j_{2}(\psi_{2})j_{1}(\psi_{1})^{-1} = \begin{pmatrix} z_{2}^{*} & -\bar{\rho} z_{1} \\ z_{1}^{*} & z_{2} \end{pmatrix}\] (here writing $z_{1}$ and $z_{2}$ for the generators of $C(S^{3}_{\theta})$), so that $eC(S^{4}_{\theta})^{4}$ is isomorphic to $M(A_{1}^{2}, A_{2}^{2}, h)$.

We pause to remark that, if $\theta=0$, then since $h^{-1}= \begin{pmatrix} z_{2} &  z_{1} \\- \rho z_{1}^{*} & z_{2}^{*} \end{pmatrix}$, the module $eC(S^{4}_{\theta})^{4}$  must be isomorphic to the module of sections of the rank-2 complex vector bundle over $S^{4}$ that is obtained by using the generator $\begin{pmatrix} z_{2} &  z_{1} \\-z_{1}^{*} & z_{2}^{*} \end{pmatrix} \in \pi_{3}(SU(2))$ as a clutching function $f_{NS}$. That vector bundle has integral of the second Chern class equal to -1 (if we use the orientation on $S^{3}$ coming from the standard orientation of $S^{4}$). Any rank-2 complex vector bundle over $S^{4}$ that supports self-dual connections with instanton number 1 is isomorphic to this bundle.

Returning to the proof that $eC(S^{4}_{\theta})^{4}$ is isomorphic to $N(2, -1)$, we note that it follows from the (proof of the) noncommutative index calculation of theorem 5.5 of \cite{natsume} that $X=\exp(2\pi it) p+1-p$ and $\begin{pmatrix} z_{1} & \sqrt{\bar{\rho}} z_{2} \\-\sqrt{\bar{\rho}} z_{2}^{*} & z_{1}^{*} \end{pmatrix}$ represent the same generator of $K_{1}(C(S^{3}_{\theta}))$. We warn the reader interested in the proof of theorem 5.5  of \cite{natsume} that, unfortunately, the matrix $Z(2)$ as given on page 1094 of \cite{natsume} is not actually in $U_{2}(C(S^{3}_{\theta}))$, as is required. Fortunately, the error is not serious and is easily fixed. The most elegant fix is to simply replace the two occurrences of $\rho_{12}$ (our $\rho$) in $Z(2)$ with $\bar\rho=\rho_{21}$, thus obtaining the above unitary matrix (There are other reasonable fixes meeting the requirements of the argument in \cite{natsume}, but all must result in unitary matrices with the same $K_{1}$-class). 

Now \[ h^{-1} = \begin{pmatrix} 1 &  0 \\ 0  & -\rho \end{pmatrix} \begin{pmatrix} 1 &  0 \\ 0 & \sqrt{\bar{\rho}} \end{pmatrix} \begin{pmatrix} z_{1} & \sqrt{\bar{\rho}} z_{2} \\-\sqrt{\bar{\rho}} z_{2}^{*} & z_{1}^{*} \end{pmatrix} \begin{pmatrix} 1 &  0 \\ 0 & \sqrt{\rho} \end{pmatrix} \begin{pmatrix} 0 &  1 \\ 1 & 0 \end{pmatrix},\] so $h^{-1}$ and $\begin{pmatrix} z_{1} & \sqrt{\bar{\rho}} z_{2} \\-\sqrt{\bar{\rho}} z_{2}^{*} & z_{1}^{*} \end{pmatrix}$ are path-connected in $GL_{2}(C(S^{3}_{\theta}))$. Thus $X$ and $h^{-1}$ are in the same $K_{1}$-class of $C(S^{3}_{\theta})$, and so $\begin{pmatrix} X &  0 \\ 0 & 1 \end{pmatrix}$ and $h^{-1}$ are path-connected in $GL_{2}(C(S^{3}_{\theta}))$, by Theorem \ref{naturalmap}.

We conclude that \[eC(S^{4}_{\theta})^{4} \cong M(A_{1}^{2}, A_{2}^{2}, h) \cong N(2,-1) \cong N(1,-1)\oplus C(S^{4}_{\theta}).\]

By similar arguments, the modules $P_{(n)}$ constructed in \cite{landivs} must be of the form $N(n+1, -(1/6)n(n+1)(n+2))$, for each $n\geq 1$. In the same way, we can show that the modules of \cite{brain} are all of the form $N(n,s)$ for $n \geq 2$.

\section{Further Directions}

I am currently attempting \cite{peterka2} to classify the finitely-generated projective modules over higher-dimensional $\theta$-deformed spheres to the extent that this is tractable. Tentatively, it appears that the situation is this: if all of the off-diagonal entries of $\theta$ are irrational, then $C(S^{l}_{\theta})$ satisfies cancellation. If, additionally, $l=2m+1$ is odd, then all $C(S^{2m-1}_{\theta})$-modules are free, and if $l=2m$ is even, then every $C(S^{2m}_{\theta})$-module is isomorphic to one of the form $N(1,s)\oplus C(S^{2m}_{\theta})^{n}$, where, as in the case $m=2$, the $C(S^{2m}_{\theta})$-module $N(1,s)$ can be thought of as a line bundle over $S^{2m}_{\theta}$ that comes from a clutching element that generates $\pi_{0}(GL_{1}(C(S^{2m-1}_{\theta}))) \cong K_{1}(C(S^{2m-1}_{\theta})) \cong \mathbb{Z}$. The generator of $\pi_{0}(GL_{1}(C(S^{2m-1}_{\theta})))$ will be one of the generators of $\pi_{0}(GL_{1}(S^{m-1}C(T^{m}_{\theta})))$, and we also will have a higher analog of Corollary \ref{tensor} in that $N(1,s)$ will be isomorphic to $M(1,s) \otimes_{C(S^{2})} C(S^{2m}_{\theta})$, where $M(1,s)$ is the module of continuous sections of a complex vector bundle over $S^{2}$ with Chern number equal to $-s$.

 But if some of the off-diagonal entries of $\theta$ are rational, then cancellation can fail. For instance, if one of the below-diagonal entries of $\theta$ is irrational, and the others are zero, then $C(S^{5}_{\theta})$ admits $\mathbb{Z}^{2}$-many ``line bundles", but no nontrivial ``higher-rank bundles", and so fails cancellation, while also differing fairly substantially from the classical case $\theta=0$. The phenomena can become quite strange for larger $l$ when some of the off-diagonal entries of $\theta$ are rational and others irrational. This contrasts with the situation for higher-dimensional noncommutative tori, in that if $\theta$ contains any irrational entries, then $C(T^{l}_{\theta})$ satisfies cancellation (see \cite{projmodules}). Of course if all of the entries of $\theta$ are rational, then $V(C(S^{l}_{\theta}))$ is isomorphic to $V(C(S^{l}))$, by Proposition \ref{moritaequiv}.

Related to the above, I have calculated some of the homotopy groups of the groups of invertible matrices over various $\theta$-deformed spheres. For example, I can show that if $\theta$ is irrational, then \[\pi_{1} (GL_{1}(C(S^{3}_{\theta}))) \cong \mathbb{Z}^{3},\] while \[\pi_{1} (GL_{k}(C(S^{3}_{\theta}))) \cong \mathbb{Z},\] for all $k\geq 2$. 

The homotopy groups of the groups of invertible matrices over $l$-dimensional noncommutative tori were calculated in \cite{homotopy} for the case that at least one entry of $\theta$ is irrational (see Theorem \ref{rieffelgenerators} above). In that case, for fixed $l$,  the homotopy groups $\pi_{j}(GL_{k}(C(T^{l}_{\theta})))$ agree for all $j\geq 0, k\geq 1$.  In general, however, the group $\pi_{j}(GL_{k}(C(S^{l}_{\theta})))$ depends partly on $\pi_{j+1}(GL_{k}(C(S^{1})))$, which in turn depends on $\pi_{j+2}(S^{m})$ for $m$ depending on $k$. Thus calculating $\pi_{j}(GL_{k}(C(S^{l}_{\theta})))$  inherits the difficulties in calculating the ordinary homotopy groups of spheres, even in the case that all off-diagonal elements of $\theta$ are irrational, and so the situation is far less well-behaved than for the noncommutative tori. But in the positive direction, the groups $\pi_{j}(GL_{k}(C(S^{l}_{\theta})))$ do stabilize somewhat faster when $\theta$ contains some irrational terms in comparison to the classical case $\theta=0$, giving hope that at least some further calculations can be made without too much difficulty.

Additionally, I am investigating the existence of instantons on the modules $N(1,s)$. These could be thought of as global abelian instantons on $S^{4}_{\theta}$. Calculations suggest that the moduli space of charge s instantons on $M(1,s)$  is a 4$|s|$-dimensional hyperk\"{a}hler space. It appears that the moduli space $M(1,s)$ has a singular component, the removal of which will leave the hyperk\"{a}hler metric incomplete, as in the classical case.

\bibliographystyle{plain}
\bibliography{arxivfinalrefs}

\end{document}